\documentclass[reqno]{amsart}
\usepackage{mathtools,amssymb,wasysym}
\usepackage{mathrsfs}
\usepackage{enumitem}
\usepackage[usenames, dvipsnames]{color}

\usepackage[letterpaper,left=1.2in,right=1.2in,top=1.5in,bottom=1.5in]{geometry}

\numberwithin{equation}{section}
\usepackage[numbers,sort&compress]{natbib}
\usepackage{hyperref}
\usepackage{esint}
\usepackage{verbatim}

\theoremstyle{plain}
\newtheorem{theorem}{Theorem}[section]
\newtheorem{lemma}[theorem]{Lemma}
\newtheorem{corollary}[theorem]{Corollary}
\newtheorem{proposition}[theorem]{Proposition}

\theoremstyle{definition}

\theoremstyle{remark}
\newtheorem{remark}[theorem]{Remark}
\newtheorem*{rem*}{Remark}

 \makeatletter
 \newcommand{\norm}{\@ifstar{\@normb}{\@normi}}
 \newcommand{\@normb}[2]{\left\Vert{#1}\right\Vert_{#2}}
 \newcommand{\@normi}[2]{\Vert{#1}\Vert_{#2}}
 \makeatother

 \global\long\def\Sob#1#2{{W}^{#1,#2}}

 \global\long\def\Leb#1{L^{#1}}

 \newcommand{\action}[1]{\left<#1 \right>}
 
 \RequirePackage{bm}

 \newcommand{\boldj}{\mathbf{j}}

 \DeclareMathOperator{\Div}{div}
 \newcommand{\relphantom}[1]{\mathrel{\phantom{#1}}}
 
 \newcommand{\myd}[1]{\,d{#1}}

 \DeclareMathOperator{\supp}{supp}

\usepackage{tcolorbox}

\makeatletter
\def\@tocline#1#2#3#4#5#6#7{\relax
  \ifnum #1>\c@tocdepth 
  \else
    \par \addpenalty\@secpenalty\addvspace{#2}%
    \begingroup \hyphenpenalty\@M
    \@ifempty{#4}{%
      \@tempdima\csname r@tocindent\number#1\endcsname\relax
    }{%
      \@tempdima#4\relax
    }%
    \parindent\z@ \leftskip#3\relax \advance\leftskip\@tempdima\relax
    \rightskip\@pnumwidth plus4em \parfillskip-\@pnumwidth
    #5\leavevmode\hskip-\@tempdima
      \ifcase #1
       \or\or \hskip 1em \or \hskip 2em \else \hskip 3em \fi%
      #6\nobreak\relax
    \hfill\hbox to\@pnumwidth{\@tocpagenum{#7}}\par
    \nobreak
    \endgroup
  \fi}
\makeatother

\allowdisplaybreaks

\begin{document}

\title[Vlasov-Riesz system]{Scattering of the Vlasov-Riesz system in the three dimensions}

\author[W. Huang]{Wenrui Huang}
\address{Department of Mathematics, Brown University, 151 Thayer Street, Providence, RI 02912, USA}
\email{wenrui\_huang@brown.edu }

\author[H. Kwon]{Hyunwoo Kwon}
\address{Division of Applied Mathematics, Brown University, 182 George Street, Providence, RI 02912, USA}
\email{hyunwoo\_kwon@brown.edu }

\thanks{\today}

\begin{abstract}
We consider an asymptotic behavior of solutions to the Vlasov-Riesz system of order $\alpha$ in $\mathbb{R}^3$ which is a kinetic model induced by Riesz interactions. We prove small data scattering when $1/2<\alpha<1$ and modified scattering when $1<\alpha<1+\delta$ for some $\delta>0$. Moreover, we show the existence of (modified) wave operators for such a regime. To the best of our knowledge, this is the first result on the existence of modified scattering with polynomial correction in kinetic models.
\end{abstract}

\maketitle

\section{Introduction}

This paper is devoted to studying the asymptotic behavior of solutions to a kinetic model induced by Riesz interactions. In this paper, we consider the following Vlasov-Riesz system of order $\alpha$: 
\begin{equation}\label{eq:VR}
\left\{
\begin{alignedat}{2}
\partial_t f + v\cdot \nabla_x f +\lambda E[f]\cdot \nabla_v f&=0&&\quad \text{in } (0,\infty)\times\mathbb{R}^3_x\times \mathbb{R}^3_v,\\
f&=f_0 &&\quad\text{on } \{t=0\}\times \mathbb{R}^3_x\times\mathbb{R}^3_v.
\end{alignedat}
\right.
\end{equation}
Here, the vector field $E[f]$ is given by
\begin{equation}\label{eq:long-range-potential}
   E[f](t,x)= \nabla_x \int_{\mathbb{R}^{3}_{y}\times\mathbb{R}^3_v} V_\alpha(x-y) f(t,y,v) \myd{y}dv, 
\end{equation}
and $\lambda\in \{-1,1\}$, where  $V_\alpha(x)$ is the \emph{Riesz kernel} defined by
\[
V_\alpha(x)=\frac{\Gamma((3-2\alpha)/2)}{\pi^{3/2}2^{2\alpha}\Gamma(\alpha)}\frac{1}{|x|^{3-2\alpha}}.
\] 
When $\alpha=1$, the kernel becomes the classical Coulomb kernel. We refer to the case $\alpha>1$ as \emph{sub-Coulombic} and the case $0<\alpha<1$ as \emph{super-Coulombic}. The case $\lambda=1$ describes an attractive case that can be observed in galactic settings. When $\lambda=-1$, \eqref{eq:VR} describes a repulsive case that can be seen in a plasma or ion gas. In this article, we call $E[f]$ the electric field for simplicity.   

We note that several researchers seek the rigorous justification of the mean limit of particle systems induced by Riesz interactions to derive fluid or kinetic equations (see e.g. Serfaty \cite{S20}, Nguyen-Rosenzweig-Serfaty \cite{NRS22}, and references therein). Riesz kernels were studied in many different physical contexts (see e.g. \cite{L22,BBDR05,DRAW02,CJ23}) and approximation theory (see e.g. \cite{HS04}). 
 
The system \eqref{eq:VR} is a natural generalization of the Vlasov-Poisson system which corresponds to $\alpha=1$. If we replace the potential $V_\alpha(x)$ with the Manev potential $\gamma |x|^{-1} +\varepsilon |x|^{-2}$
then the system becomes the Vlasov-Manev system. The Manev potential was introduced by Manev \cite{M24,M25,M30,M30b} to explain the perihelion of Mercury without using Einstein's theory of general relativity. When $\gamma=0$ and $\varepsilon \neq 0$, the system corresponds to the Vlasov-Riesz system of order $1/2$. For simplicity, we will call this system the \emph{pure Vlasov-Manev system}.

For the case of the Vlasov-Poisson system ($\alpha=1$), Bardos-Degond \cite{BD85} proved the global existence of smooth solutions with small initial data \cite{BD85}. Later, Lions-Perthame \cite{LP91} and Pfaffelmoser \cite{P92} independently proved the global existence of classical solutions for large initial data. Hence, it is natural to ask for the asymptotic behavior of solutions; in particular, we are interested in whether the solution scatters.  In this direction, Choi-Ha \cite{CH11} proved that any nontrivial classical solution of the Vlasov-Poisson equation does not have a linear scattering, that is, a nonlinear solution cannot converge to a solution of the linear equation as time goes to infinity. On the other hand, Choi-Kwon \cite{CK16} proved the existence of the modified scattering of the Vlasov-Possion system, that is, the solution to a nonlinear equation converges to a solution of a linear equation with a suitable logarithmic correction. However, in their result, the trajectory was not explicit. 

Recently, there have been several works on kinetic equations which employ  techniques from the analysis of dispersive equations. For relevant results to our problem, Ionescu-Pausader-Wang-Widmayer \cite{IPWW22} proved modified scattering of solutions with a concrete scattering trajectory by using so-called $Z$-norm method, which was used in several papers to show the global stability of various models. Later, Flynn-Ouyang-Pausader-Widmayer \cite{FOPW23} proved a similar result by using pseudo-conformal transformation motivated by \cite{CN18, NAYO1987,T09} (see also references therein). Moreover, they proved the existence of wave operators and scattering maps. {{We also refer to \cite{BiVe} and \cite{ScTa} for more precise expansions of the asymptotic behavior and to \cite{HuPaSu} for the asymptotic behavior in the presence of a boundary.}}

 For the case of the Vlasov-Manev system ($\alpha=1/2$), Bobylev et al. \cite{BDIV97,BDIV98} proved the local existence of classical solutions and finite-time singularity formation of smooth solutions for the gravitational case. Later, for the pure Vlasov-Manev system, Lemou-M\'ehats-Rigault \cite{LMR12} proved the existence of stable ground states, and as an application, they proved the existence of spherically symmetric self-similar blow-up for the gravitational case. 

However, there are few results concerning the Vlasov-Riesz system. Recently, Choi-Jeong \cite{CJ23}  proved the {local-wellposedness} of \eqref{eq:VR} in a weighted Sobolev space on $\mathbb{R}^3$ when $\alpha>3/8$. Moreover, for the gravitational case, they proved a finite-time singularity formation of solutions when $3/8<\alpha\leq 1/2$. They also proved the same result by replacing the vector field $E[f]$ with a more general form, which extends the result of Bobylev et al. \cite{BDIV98}. { In a different paper, Choi-Jeong \cite{CJ23b} constructed global weak solutions to the Vlasov-Manev-Fokker-Planck system which also covers Vlasov-Riesz system. However, it is difficult to analyze the asymptotic behavior of the solution in their framework. }Hence, it is natural to ask about the range on $\alpha$ which guarantees a global existence of strong solutions to \eqref{eq:VR} that scatter. Once we have scattering results, it is reasonable to investigate the existence of a (modified) wave operator.

The purpose of this paper is to give an affirmative answer to these questions. We show that if $\alpha \in (1/2,1+\delta)\setminus\{1\}$ for some $\delta>0$, then the solution exists globally {for small initial data}. For asymptotic behavior of the global solution, for the super-Coulombic case, if $1/2<\alpha<1$, then the solution has a linear scattering. For the sub-Coulombic case, $1<\alpha<1+\delta$, the solution has a modified scattering with polynomial corrections. Moreover, we prove the existence of wave operator when $1/2<\alpha<1$ and modified wave operator when $1<\alpha<1+\delta$. To the best of our knowledge, our result is the first example of solutions to kinetic models that exhibit modified scattering with polynomial corrections.

To elaborate the main result of this paper, we reformulate our problem \eqref{eq:VR} into \eqref{eq:VR-mu} by writing $\mu^2=f$ as in \cite{IPWW22,FOPW23}:
\begin{equation}\label{eq:VR-mu}
\left\{
\begin{alignedat}{2}
\partial_t \mu + v\cdot \nabla_x \mu + \lambda E[\mu]\cdot \nabla_v \mu&=0&&\quad \text{in } (0,\infty)\times\mathbb{R}^3_x\times \mathbb{R}^3_v,\\
\mu&=\mu_0 &&\quad\text{on } \{t=0\}\times \mathbb{R}^3_x\times\mathbb{R}^3_v.
\end{alignedat}
\right.
\end{equation}
where 
\begin{equation*}
	E[\mu](t,x) = \nabla_x \iint_{\mathbb{R}^3_y\times\mathbb{R}^3_v} V_\alpha(x-y) \mu^2(t,y,v)dvdy.
\end{equation*}

Now, we are ready to state the main results of this paper.  Our first result is the existence of modified scattering of the problem \eqref{eq:VR-mu}.
\begin{theorem}\label{thm:A}
There exists a constant $\delta>0$ so that the following is true for $\alpha \in (1/2,1+\delta)\setminus\{1\}$: {there exists a small $\varepsilon>0$ such that for any $\mu_0$ satisfying
\begin{equation*}\label{eq:mu1-control}
   \norm{\mu_0}{\Leb{2}_{x,v}}+\norm{\action{x}^4\mu_0}{\Leb{\infty}_{x,v}}+\norm{\action{v}^4\mu_0}{\Leb{\infty}_{x,v}}+\norm{\nabla_{x,v}\mu_0}{\Leb{2}_{x,v}\cap\Leb{\infty}_{x,v}}\leq\varepsilon,\end{equation*}
there exists a unique global strong solution $\mu$ of the initial value problem for \eqref{eq:VR-mu} with $\mu(0,x,v)=\mu_0(x,v)$. Moreover, there exist $\mu_\infty\in  L^2_{x,v}\cap\Leb{\infty}_{x,v}$ and $E_\infty \in \Leb{\infty}_v$}  such that {locally} uniformly in $(x,v)$,
\begin{equation}\label{eq:VR-asymptotic}
  \mu\left(t,x+tv+\frac{\lambda}{2-2\alpha}t^{2\alpha-2} E_\infty(v),v\right)\rightarrow \mu_\infty(x,v)\quad \text{as } t\rightarrow\infty.
  \end{equation}
\end{theorem}

\begin{remark}\leavevmode 
\begin{enumerate} 
\item For the case of the Vlasov-Poisson system ($\alpha=1$), Flynn-Ouyang-Pausader-Widmayer \cite{FOPW23} proved that the solution of the Vlasov-Poisson system has the modified scattering dynamics with logarithmic corrections:
\[  
\mu\left(t,x+tv-\lambda \ln(t) E_\infty(v),v\right)\rightarrow \mu_\infty(x,v)\quad \text{as } t\rightarrow \infty.\]
Compared to this result, our result shows the existence of a linear scattering when $1/2<\alpha<1$ and modified scattering with polynomial correction when $\alpha>1$. 
\item When {$\alpha>7/8$}, the term {$\norm{\action{x}^4\mu_0}{\Leb{\infty}_{x,v}}$ and $\norm{\action{v}^4\mu_0}{\Leb{\infty}_{x,v}}$} in \eqref{eq:VR-mu} can be replaced by {$\norm{\action{x}^2\mu_0}{\Leb{\infty}_{x,v}}$ and $\norm{\action{v}^2\mu_0}{\Leb{\infty}_{x,v}}$}. It is necessary to control higher weight on the initial data when we want to obtain a scattering if {$1/2<\alpha \leq 7/8$.} 
\item Our theorem requires us to assume that $\alpha \in (1/2,1+\delta)\setminus\{1\}$ for small $\delta>0$. In fact, we can show that global solution exists for $\alpha \in (1/2,10/9)\setminus\{1\}$ and with the help of Mathematica, the solution has a modified scattering for $\alpha \in (1,1.08]$.  Hence one might ask whether we can extend the result for all $\alpha \in (1/2,3/2)\setminus\{1\}$, see Remark \ref{remark:threshold} for further discussion.

In the dispersive settings, Ginibre-Velo \cite{GV00a,GV00b} and Hayashi-Naumkin \cite{HN01} obtained modified scattering with polynomial correction for modified Hartree-type equations and Hartree-type equations, respectively.
\end{enumerate} 
\end{remark}

The next theorem concerns the existence of (modified) wave operators and scattering maps for the Vlasov-Riesz system. To describe this, we define 
\begin{equation}\label{eq:Asymptoric-electric-field}
	E_\infty[\mu](v):=\nabla_v \iint_{\mathbb{R}^3_y\times\mathbb{R}^3_w} V_{\alpha}(v-w) \cdot \mu^2(y,w)\myd{y}dw.
\end{equation}

\begin{theorem}\label{thm:B}
There exists a constant $\delta>0$ such that for $\alpha \in (1/2,1+\delta)\setminus\{1\}$, the following hold: 
\begin{enumerate}[label=\textnormal{(\roman*)}]
\item {Let $\mu_\infty \in \Sob{2}{\infty}(\mathbb{R}^3_x\times \mathbb{R}^3_v)$ and $E_\infty=E_\infty	[\mu_{\infty}]  \in \Sob{3}{\infty}(\mathbb{R}^3_v)$ satisfy 
\begin{equation*}\label{eq:assumption-wave}
\begin{aligned}
A&=\norm{\mu_\infty}{\Leb{2}_{x,v}}+\norm{\action{x}^5\mu_\infty}{\Leb{2}_{x,v}\cap\Leb{\infty}_{x,v}}+\norm{\action{x}\nabla_{x,v} \mu_\infty}{\Leb{2}_{x,v}\cap\Leb{\infty}_{x,v}}\\
&\relphantom{=}+\norm{\action{x}^2\nabla^2_{x,v}\mu_\infty}{\Leb{\infty}_{x,v}}+\norm{E_\infty}{W^{3,\infty}_{x,v}}<\infty.
\end{aligned}
\end{equation*}
\begin{enumerate}
\item If $\alpha>1$, there exists a global strong solution $\mu$ of \eqref{eq:VR} satisfying \eqref{eq:VR-asymptotic}. \\
\item there exists a small $\varepsilon_0>0$ such that if $A\leq \varepsilon_0$, then there exists a unique global strong solution $\mu$ of \eqref{eq:VR} satisfying \eqref{eq:VR-asymptotic}. 
\end{enumerate}
}
\item {there exists a small $\varepsilon_0>0$} such that for any $\mu_{-\infty} \in \Sob{2}{\infty}(\mathbb{R}^3_x\times\mathbb{R}^3_v)$ and $E_{-\infty}=E_{\infty}[\mu_{-\infty}]  \in \Sob{3}{\infty}(\mathbb{R}^3_v)$ satisfying
\begin{equation*}\label{eq:assumption-scattering-map}
\norm{\mu_{-\infty}}{\Leb{2}_{x,v}}+\norm{\action{x,v}^5\mu_{-\infty}}{L^2_{x,v}\cap\Leb{\infty}_{x,v}}+\norm{\action{x}\nabla_{x,v}\mu_{-\infty}}{L^2_{x,v}\cap\Leb{\infty}_{x,v} }+\norm{\action{x}^2\nabla^2_{x,v}\mu_{-\infty}}{\Leb{\infty}_{x,v}}\leq \varepsilon,
\end{equation*}
there exist a unique global strong solution $\mu$ of \eqref{eq:VR}, $\mu_\infty \in \Leb{2}_{x,v} \cap \Leb{\infty}_{x,v}$ and {$E_\infty(v)\in \Leb{\infty}_v $} such that 
\begin{equation*}
\mu\left(t,x+tv\pm \frac{\lambda}{2-2\alpha} t^{2\alpha-2} E_{\pm \infty}(v),v\right)\rightarrow \mu_{\pm \infty}(x,v)\quad \text{as } t\rightarrow\pm \infty.
\end{equation*}
\end{enumerate} 
\end{theorem}

\begin{remark}\leavevmode
\begin{enumerate}
\item When $\alpha=1$, Flynn-Ouyang-Pausader-Widmayer \cite{FOPW23} proved the existence of modified wave operator and scattering map.
\item A more precise statement is given in Theorem \ref{thm:Wave}. The restriction on the range $\alpha$ is due to the convergence rate of the electric field (see Proposition \ref{prop:continuity-of-E-2}). {We do not need to assume smallness on $A$ when $\alpha>1$. However, it is unclear whether the solution is unique without smallness. It is interesting to consider whether we can construct a wave operator when $\alpha<1$ without the smallness assumption; see Remark \ref{rem:wave-operator-large}.}
\item As an analogy between Hartree equations and Vlasov-Riesz system, Ginibre-Ozawa \cite{GO93} first obtained existence of modified wave operator for Hartree equations when the potential is Coulombic. The general case was obtained by Nawa-Ozawa \cite{NO92} and Nakanishi \cite{N02A,N02B} for a certain range.
\end{enumerate} 
\end{remark}

Let us outline the proof of Theorems \ref{thm:A} and \ref{thm:B}. Inspired by an analogy between dispersive equations and kinetic equations as in \cite{IPWW22,FOPW23,LMR12}, we use a pseudo-conformal symmetry of \eqref{eq:VR} which could make it easy to analyze global dynamics of solutions to \eqref{eq:VR} by considering the problem on $(0,1]$. We prove the local-in-time solution becomes global via a standard bootstrap argument by considering the equation as a nonlinear transport equation. {Unlike \cite{FOPW23}, the Vlasov-Riesz system has different features compared to the Vlasov-Poisson system.  In the classical Vlasov-Poisson system, for instance, the moments grow at most logarithmically and the gradient of the electric field decays like {$\action{t}^{-3}|\ln t|^p$}, allowing for logarithmic errors to be controlled easily.} When $\alpha<1$,  since the interaction kernel has a shorter range effect than that of the Vlasov-Poisson system, it is reasonable to get the linear scattering result as in the Hartree equations with a short-range interaction kernel \cite{NAYO1987}. However, we need to impose more weight on $\mu$ to guarantee the existence of scattering. Our main novelty lies in the case $\alpha>1$, where a delicate analysis is required to establish global existence and modified scattering. This difficulty arises because the $L^\infty-$weighted norm of the solution exhibits polynomial growth in time, while at the same time the gradient of the electric field decays at a slower rate. This competition forces us to obtain global existence and its modified scattering in a restricted regime $1<\alpha<1+\delta$ for small $\delta>0$ { to close the bootstrap argument.}

The proof of Theorem \ref{thm:B} also uses pseudo-conformal symmetry of the system. The transformation converts the original problem into the local well-posedness problem starting from $s=0$. However, due to the pseudo-conformal inversion, the problem is not well-defined at $s=0$. To fix this issue, we find a change of variables for spatial and velocity variables, which preserves the symplectic structures so that it preserves the Hamiltonian structure of the Vlasov-Riesz system and the new Hamiltonian {is} less singular at $s=0$. This enables us to consider the problem at $s=0$ to show the existence of the wave operator. To show the local well-posedness of the transformed problem, we obtain a priori estimates for solutions that involve momentum estimates and derivative estimates. Unlike the original local well-posedness problem (see Theorem \ref{thm:VR-LWP}), we need to mitigate the singularity at $s=0$ when we estimate derivatives of the transformed solutions, which makes us introduce appropriate weights to estimate those. Then the solution can be obtained via a standard Picard iteration (Theorem \ref{thm:Wave}). Once we obtain a local solution, using the pseudo-conformal transformation again, the local-in-time solution becomes a strong solution to \eqref{eq:VR} on $[T,\infty)$ for some $T>0$. {When $\alpha>1$, we apply Theorem \ref{thm:propagation} to construct a modified wave operator using time-reversal symmetry. If we assume smallness, then we apply Theorem \ref{thm:A} to construct wave operators.} Finally, Theorem \ref{thm:B} (ii) follows by combining Theorems \ref{thm:A} { and \ref{thm:B} (i)(b).}

The organization of this paper is as follows. Section \ref{sec:2} contains several facts on pseudo-conformal transformation, transport equations, and control of the electrical field $E$. These controls play a crucial role in the proof of the global existence of solutions as well as the existence of (modified) wave operators. Theorem \ref{thm:A} will be proved in Section \ref{sec:3}. In Section \ref{sec:4}, we introduce a change of variable that will be used to construct a wave operator and show the local well-posedness result that is crucial to constructing a wave operator and scattering map for \eqref{eq:VR}. For the sake of completeness, we include a proof of the local well-posedness of the Vlasov-Riesz system in the Appendix \ref{app:A} that we need.  We also provide a proof of interpolation inequalities in Appendix \ref{app:interpolation} that will be used in the paper. 

Finally, we introduce some notation used in this paper. Let $\mathbb{R}^3$ denote the standard Euclidean space of points $x=(x^1,x^2,x^3)$ and let $\action{x}=(1+|x|^2)^{1/2}$. We frequently suppress the time variable when it involves norms, e.g.,  
\[ \norm{\mu(s)}{\Leb{r}_{q,p}}=\norm{\mu}{\Leb{r}_{q,p}}.\]
For $R>0$, let $B_R$ be the Euclidean open ball with radius $R$ centered at the origin. {  For two nonnegative quantities $A$ and $B$, we write $A\apprle_{\alpha,\beta,\dots} B$ if $A\leq CB$ for some positive constant $C$ depending on the parameters $\alpha$, $\beta$, ... If the dependence is evident, we usually omit the subscripts and simply write $A\apprle B$. We write $A\approx B$ if $A\apprle B$ and $B\apprle A$.  For convenience, we write $A\ll 1$ if $A$ is sufficiently small. Throughout this paper, the time variable $s$ will be always less than $2$ and  we will frequently use the fact that for any constant  $\kappa>0$, $\max\{1, s^{-\kappa}\}\approx s^{-\kappa}$ and $\max\{1,s^{\kappa}\} \approx 1$.}

\subsection*{Acknowledgement}
H. Kwon thanks  Prof. Yonggeun Cho for a detailed explanation on modified scattering results in dispersive equations by introducing \cite{GV00a,GV00b} to him. The authors thank Prof. Beno\^{i}t Pausader for several comments and encouragement while we are preparing this paper. Finally, the authors would like to express sincere gratitude to anonymous referees for valuable comments and suggestions, which have greatly improved the manuscript.

\section{Preliminaries}\label{sec:2}
In this section, we first introduce the pseudo-conformal transformation related to the Vlasov-Riesz system and introduce the transport equation associated with the Vlasov-Riesz system. Finally, we give some estimates on the electric field associated with the solution.

\subsection{Pseudo-conformal transformation}
Recall the involution of $\mathbb{R}\setminus\{0\}\times\mathbb{R}^3\times\mathbb{R}^3$ given by the \emph{pseudo-conformal inversion} 
\[  \mathcal{I}:(t,x,v)\mapsto \left(\frac{1}{t},\frac{x}{t},x-tv\right).\]

Define 
\begin{equation}\label{eq:pseudo-conformal-relationship}
 \gamma(s,q,p):=\mu\left(\frac{1}{s},\frac{q}{s},q-sp\right)\quad\text{and}\quad \mu(t,x,v)=\gamma\left(\frac{1}{t},\frac{x}{t},x-tv\right).
\end{equation}
Then 
\[  E[\mu](t,tx)=\frac{1}{t^{4-2\alpha}}E[\gamma]\left(\frac{1}{t},x\right).\]
and it is easy to show that $\mu$ solves \eqref{eq:VR-mu} on $0\leq T_*\leq t\leq T^*$ if and only if $\gamma$ satisfies 
\begin{equation}\label{eq:VR-pct}
  (\partial_s + p\cdot\nabla_q)\gamma + \lambda s^{1-2\alpha}  E[\gamma]\cdot\nabla_p \gamma=0
\end{equation}
on $0\leq (T^*)^{-1} \leq s\leq (T_*)^{-1}$.

\subsection{Transport equations}
We can rewrite \eqref{eq:VR-pct} into the transport equation
\begin{equation*}\label{eq:VR-transport}
\mathcal{L}[\gamma]=\partial_s \gamma + \Div_{q,p} \{(p,\lambda s^{1-2\alpha} E(q))\cdot \gamma\}=0.
\end{equation*}
Hence if $h$ is a strong solution of 
\begin{equation*}
\mathcal{L}[h]=F(s,q,p)
\end{equation*}
with $h(1) \in \Leb{r}_{q,p}$ for some $r\geq 1$, then 
\begin{equation}\label{eq:estimate-transport}
\norm{h(s)}{\Leb{r}_{q,p}}\leq \norm{h(1)}{\Leb{r}_{q,p}}+\int_s^1 \norm{F(s')}{\Leb{r}_{q,p}}ds'
\end{equation}
for some $0\leq s\leq 1$ in the interval of existence. 

Define 
\[ \mathcal{L}_{\gamma}[h]=\partial_s h + p\cdot \nabla_q h + \lambda s^{1-2\alpha} E[\gamma]\cdot \nabla_p h.
\]
If $\gamma$ is a solution to \eqref{eq:VR-pct}, then for any $m,n\in \{1,2,3\}$, we have 
\begin{equation}\label{eq:commutator-relationship}
\begin{gathered}
\mathcal{L}_{\gamma}[q^m \gamma]=\mathcal{L}_{\gamma}[q^m]\gamma = p^m \gamma,\quad \mathcal{L}_{\gamma}[p^m\gamma]=\lambda s^{1-2\alpha} E^m[\gamma] \gamma,\\
\mathcal{L}_{\gamma}[\partial_{q^m} \gamma] = \partial_{q^m} (\mathcal{L}_{\gamma}[\gamma])-(\partial_{q^m} \mathcal{L}_{\gamma})[\gamma]=-\lambda s^{1-2\alpha} \partial_{q^m} E^j[\gamma] \partial_{p^j} \gamma,\\
\mathcal{L}_{\gamma}[\partial_{p^m} \gamma]=-\partial_{q^m} \gamma.
\end{gathered}
\end{equation}

\subsection{Controls on electric fields}
{ Let $\chi \in C_c^\infty(B_2)$ be a radially symmeric function such that $\supp \chi \subset B_2\setminus B_{1/2}$ and $\int_{\mathbb{R}^3} \chi \myd{x}=1$. By a change of variable, we have
\[  \frac{1}{|x-y|^{3-2\alpha}} = c\int_0^\infty R^{-(3-2\alpha)} \chi(R^{-1}(x-y))\frac{dR}{R}\]
for some constant $c>0$.
Then we can rewrite $E$ into
\begin{equation}\label{eq:integral-expression}
\begin{aligned}
E[\gamma](s,q) &=c\int_0^\infty \left[\iint_{\mathbb{R}^3_{y}\times\mathbb{R}^3_p} R^{-(3-2\alpha)}(\nabla \chi)(R^{-1}(q-y))\gamma^2(s,y,p)dydp\right]\frac{dR}{R^2}\\
&=\int_0^\infty E_R(s,q)\frac{dR}{R^2},
\end{aligned}
\end{equation}
 where 
\[  E_R(s,q) = c_1\iint_{\mathbb{R}^3_y \times \mathbb{R}^3_p}   R^{-(3-2\alpha)}(\nabla \chi)(R^{-1}(q-y))\gamma^2(s,y,p)dpdy\]
for some constant $c_1>0$.

One can easily see that for $k=0,1,2$, we have
\begin{equation}\label{eq:ER-estimate}
 |\nabla^k_q E_R(s,q)|\apprle R^{-(3-2\alpha)-k}\norm{\gamma}{\Leb{2}_{q,p}}^2,
\end{equation}
which are useful for large $R$.}

To estimate the integral \eqref{eq:integral-expression} for small $R$, we decompose $E_R$ into 
\[   E_{R}(s,q)=\int_0^\infty E_{R,V}(s,q) \frac{dV}{V},\]
{where 
\[   E_{R,V}(s,q)=c_1\iint_{\mathbb{R}^3_y\times \mathbb{R}^3_u} R^{-(3-2\alpha)} (\nabla \chi)(R^{-1}(q-y)) \chi(V^{-1} u) \gamma^2(s,y,u) dudy.\]
Then note that }
\begin{equation}\label{eq:ERV-estimate}
\begin{aligned}
|E_{R,V}(s,q)|&\apprle_\alpha R^{2\alpha} \min\{ V^3 \norm{\gamma}{\Leb{\infty}_{q,p}}^2, V^{-1-2N}\norm{|p|^{2+N} \gamma}{\Leb{\infty}_{q,p}}^2, R{V^{-1}} \norm{|p|^4 \gamma}{\Leb{\infty}_{q,p}}\norm{\nabla_q \gamma}{\Leb{\infty}_{q,p}} \},\\
|\nabla E_{R,V}(s,q)|&\apprle_\alpha R^{2\alpha-1} \min\left\{RV^3 \norm{\nabla_q \gamma}{\Leb{\infty}_{q,p}}\norm{\gamma}{\Leb{\infty}_{q,p}},RV^{-1} \norm{|p|^4 \gamma}{\Leb{\infty}_{q,p}}\norm{\nabla_q \gamma}{\Leb{\infty}_{q,p}},V^{-1}\norm{|p|^2\gamma}{\Leb{\infty}_{q,p}}^2\right\}.
\end{aligned}
\end{equation}
for any $R,V>0$ and {for all} $N\in \mathbb{N}\cup\{0\}$. { Indeed, a change of variable gives
\begin{align*}
|E_{R,V}(s,q)|&\apprle \norm{\gamma}{\Leb{\infty}_{q,p}}^2\iint_{\mathbb{R}^3_y\times\mathbb{R}^3_u} R^{-(3-2\alpha)} |(\nabla \chi)(R^{-1}(q-y))|\chi(V^{-1}u)\myd{u}dy\\
&\apprle_\alpha\norm{\gamma}{\Leb{\infty}_{q,p}}^2 R^{2\alpha} V^3.
\end{align*}
for all $(s,q)\in I\times \mathbb{R}^3$. For the second inequality, for any $N\in\mathbb{N}\cup\{0\}$, a change of variable gives
\begin{align*}
|E_{R,V}(s,q)|&\apprle \norm{|p|^{2+N}\gamma}{\Leb{\infty}_{q,p}}^2\iint_{\mathbb{R}^3_y\times\mathbb{R}^3_u} R^{-(3-2\alpha)} |(\nabla \chi)(R^{-1}(q-y))|\chi(V^{-1}u)|u|^{-4-2N}\myd{u}dy\\
&\apprle_\alpha \norm{|p|^{2+N}\gamma}{\Leb{\infty}_{q,p}}^2 R^{-(3-2\alpha)}  R^3 V^{-1-2N}.
\end{align*}
for all $(s,q)\in I\times \mathbb{R}^3$. For the third inequality, integration by part gives
\begin{align*}
E_{R,V}(s,q)&=-2c_1\iint_{\mathbb{R}^3_y\times\mathbb{R}^3_u} R^{-(2-2\alpha)} \chi(R^{-1}(q-y)) \chi(V^{-1}u) \nabla_y \gamma(s,y,u). \gamma(s,y,u)dudy.
\end{align*}
Then a change of variable gives
\begin{align*}
|E_{R,V}(s,q)|&\apprle \norm{|p|^4 \gamma}{\Leb{\infty}_{q,p}}\norm{\nabla_q \gamma}{\Leb{\infty}_{q,p}}\int_{\mathbb{R}^3_{y}\times\mathbb{R}^3_u} R^{-(2-2\alpha)} \chi(R^{-1}(q-y)) \chi(V^{-1}u) |u|^{-4}dudy\\
&\apprle R^{2\alpha+1} V^{-1}\norm{|p|^4 \gamma}{\Leb{\infty}_{q,p}}\norm{\nabla_q \gamma}{\Leb{\infty}_{q,p}}.
\end{align*}
One can prove the desired inequality for $\nabla E_{R,V}$ in a similar way.
}

{{As a corollary, using the first two estimates of the first line in \eqref{eq:ERV-estimate}, we have 
\begin{equation}\label{eq:ER-estimate-3}
	|E_R(s,q)|\lesssim R^{2\alpha}\int_0^1 V^3 \norm{\gamma}{L^\infty_{q,p}}^2 \frac{dV}{V}+R^{2\alpha}\int_1^\infty V^{-1}\norm{|p|^2\gamma}{L^\infty_{q,p}}^2\frac{dV}{V}\lesssim R^{2\alpha}\norm{\action{p}^2\gamma}{L^{\infty}_{q,p}}^2.
\end{equation}}}

The following proposition shows that under suitable assumption on $\alpha$ and the regularity of functions $\gamma$, we can show that the electric field and its gradient are bounded.

\begin{proposition}\label{prop:electric-field-boundedness}
Let $1/2<\alpha <3/2$ and $A>0$.
\begin{enumerate}[label=\textnormal{(\roman*)}]
\item if $\gamma \in \Leb{\infty}_s\Leb{2}_{q,p}$ and $|p|^2\gamma \in \Leb{\infty}_{s,q,p}$, then
\begin{align*}
|E(s,q)|&\apprle A^{4-2\alpha} \left[ \norm{\gamma}{\Leb{2}_{q,p}}^2 + \norm{\gamma}{\Leb{\infty}_{q,p}}^2\right] + A^{-2\alpha} \norm{|p|^2 \gamma}{\Leb{\infty}_{q,p}}^2,
\end{align*} 
and {
\begin{align*}
|E(s,q)|&\apprle A^{4-2\alpha}\left[\norm{\gamma}{\Leb{2}_{q,p}}^2 + \norm{\gamma}{\Leb{\infty}_{q,p}}^2\right] + A^{-2\alpha-1} \norm{|p|^4\gamma}{\Leb{\infty}_{q,p}}\norm{\nabla_q \gamma}{\Leb{\infty}_{q,p}}.
\end{align*} }
\item if $1-2\alpha<\theta<(2\alpha-1)/3$ and if $\gamma \in \Leb{\infty}_{s,p}\Sob{1}{\infty}_{q}$, and $|p|^4 \gamma \in \Leb{\infty}_{s,q,p}$ in addition, then {
\begin{equation}\label{eq:gradient-E}
\begin{aligned}
|\nabla_q E(s,q)|&\apprle A^{5-2\alpha}\norm{\gamma}{\Leb{2}_{q,p}}^2 + A^{1-2\alpha+3\theta}\norm{\gamma}{\Leb{\infty}_{q,p}}\norm{\nabla_q \gamma}{\Leb{\infty}_{q,p}}\\
&\relphantom{=} + A^{1-2\alpha-\theta} \norm{|p|^4\gamma}{\Leb{\infty}_{q,p}}\norm{\nabla_q \gamma}{\Leb{\infty}_{q,p}}.
\end{aligned}
\end{equation} 
\item if $2-2\alpha<\theta<(2\alpha-1)/3$ and if $\gamma \in \Leb{\infty}_{s,p}\Sob{1}{\infty}_{q}$, and $|p|^2 \gamma \in \Leb{\infty}_{s,q,p}$ in addition, then 
\begin{equation}\label{eq:gradient-E-2}
\begin{aligned}
|\nabla_q E(s,q)|&\apprle A^{5-2\alpha}\norm{\gamma}{\Leb{2}_{q,p}}^2 + A^{1-2\alpha+3\theta}\norm{\gamma}{\Leb{\infty}_{q,p}}\norm{\nabla_q \gamma}{\Leb{\infty}_{q,p}}+ A^{2-2\alpha-\theta} \norm{|p|^2\gamma}{\Leb{\infty}_{q,p}}^2.
\end{aligned}
\end{equation} }
\end{enumerate}
\end{proposition}
\begin{remark}
If we choose $A=1$ in \eqref{eq:gradient-E}, then 
\begin{equation*}\label{eq:gradient-control-less-than-1}
 \norm{\nabla_q E}{\Leb{\infty}}\apprle \norm{\gamma}{\Leb{2}_{q,p}}^2 + \norm{\gamma}{\Leb{\infty}_{q,p}}\norm{\nabla_q \gamma}{\Leb{\infty}_{q,p}} + \norm{|p|^4 \gamma}{\Leb{\infty}_{q,p}}\norm{\nabla_q \gamma}{\Leb{\infty}_{q,p}}.
\end{equation*}
This will be used for the case $1/2<\alpha \leq 7/8$. Similarly, it follows from \eqref{eq:gradient-E-2} that for $\beta<0$, we have
\begin{equation}\label{eq:gradient-control-less-than-uniform}
 \norm{\nabla_q E}{\Leb{\infty}}\apprle \norm{\gamma}{\Leb{2}_{q,p}}^2 +\norm{\gamma}{\Leb{\infty}_{q,p}}\norm{\nabla_q \gamma}{\Leb{\infty}_{q,p}} + \norm{|p|^2 \gamma}{\Leb{\infty}_{q,p}}^2
\end{equation}
and
\begin{equation}\label{eq:gradient-control-less-than-2}
\begin{aligned}
 \norm{\nabla_q E}{\Leb{\infty}}&\apprle s^{\beta(5-2\alpha)}\norm{\gamma}{\Leb{2}_{q,p}}^2 +s^{\beta(1-2\alpha+3\theta)} \norm{\gamma}{\Leb{\infty}_{q,p}}\norm{\nabla_q \gamma}{\Leb{\infty}_{q,p}} +{ s^{\beta(2-2\alpha-\theta)}\norm{|p|^2 \gamma}{\Leb{\infty}_{q,p}}^2.}
\end{aligned}
\end{equation}
We will use \eqref{eq:gradient-control-less-than-uniform} when $7/8<\alpha<1$  and \eqref{eq:gradient-control-less-than-2} when $\alpha>1$.
\end{remark}

\begin{proof}
(i) We decompose it into  
\begin{equation}\label{eq:E-representation-V-R}
  E(s,q)=\int_0^A\int_0^B E_{R,V}\frac{dV}{V} \frac{dR}{R^2}+\int_0^A\int_B^\infty E_{R,V}\frac{dV}{V} \frac{dR}{R^2} + \int_A^\infty E_R \frac{dR}{R^2}.
  \end{equation}
Then by \eqref{eq:ER-estimate} and \eqref{eq:ERV-estimate}, we have {
\begin{equation}\label{eq:E-estimate}
|E(s,q)|\apprle_\alpha A^{2\alpha-1} B^3 \norm{\gamma}{\Leb{\infty}_{q,p}}^2+ A^{2\alpha-1}B^{-1} \norm{|p|^2\gamma}{\Leb{\infty}_{q,p}}^2+A^{2\alpha-4} \norm{\gamma}{\Leb{2}_{q,p}}^2.
\end{equation}
Hence the desired estimate follows by choosing $A^{-1}$ and $A$ instead of $A$ and $B$, respectively. }

To show the second assertion, we use the same representation \eqref{eq:E-representation-V-R} and inequalities  \eqref{eq:ER-estimate} and \eqref{eq:ERV-estimate} to get
\[ |E(s,q)|\apprle_\alpha A^{2\alpha-1}B^3\norm{\gamma}{\Leb{\infty}_{q,p}}^2+A^{2\alpha}B^{-1}\norm{|p|^4\gamma}{\Leb{\infty}_{q,p}}\norm{\nabla_q \gamma}{\Leb{\infty}_{q,p}}+A^{2\alpha-4}\norm{\gamma}{\Leb{2}_{q,p}}^2. \]
By choosing $A^{-1}$ and $A$ instead of $A$ and $B$, we get the desired result.

(ii) We write 
\begin{equation}
\begin{aligned}
\nabla E(s,q)&=\int_0^A \int_0^{R^{-\theta}} \nabla E_{R,V}(s,q)\frac{dV}{V}\frac{dR}{R^2} +\int_0^A \int_{R^{-\theta}}^\infty \nabla E_{R,V}(s,q)\frac{dV}{V}\frac{dR}{R^2} \\
&\relphantom{=}+ \int_A^\infty \nabla E_{R}(s,q)\frac{dR}{R^2}.
\end{aligned}
\end{equation}
Since $1-2\alpha<\theta<(2\alpha-1)/3$, { it follows from \eqref{eq:ER-estimate} and \eqref{eq:ERV-estimate} that
\begin{equation}\label{eq:gradient-E-proof}
\begin{aligned}
|\nabla E(s,q)|&\apprle_\alpha  A^{2\alpha-3\theta-1} \norm{\nabla_q \gamma}{\Leb{\infty}_{q,p}} \norm{ \gamma}{\Leb{\infty}_{q,p}}+A^{2\alpha-1+\theta}\norm{|p|^4\gamma}{\Leb{\infty}_{q,p}}\norm{\nabla_q\gamma}{\Leb{\infty}_{q,p}}+ A^{2\alpha-5} \norm{\gamma}{\Leb{2}_{q,p}}^2.
\end{aligned}
\end{equation}
Hence the desired estimate follows by choosing $A^{-1}$ instead of $A$.

(iii) Using $R^{2\alpha-1}V^{-1}\norm{|p|^2 \gamma}{\Leb{\infty}_{q,p}}^2$ instead, the proof is similar except calculating the second integral in \eqref{eq:gradient-E-proof} differently. Then we get 
\[
 |\nabla E(s,q)|\apprle_\alpha A^{2\alpha-3\theta-1} \norm{\nabla_q\gamma}{\Leb{\infty}_{q,p}}\norm{\gamma}{\Leb{\infty}_{q,p}}+A^{2\alpha-2+\theta}\norm{|p|^2\gamma}{\Leb{\infty}_{q,p}}^2+A^{2\alpha-5}\norm{\gamma}{\Leb{2}_{q,p}}^2.
\]
Hence the desired estimate follows by choosing $A^{-1}$ instead of $A$. This completes the proof of Proposition \ref{prop:electric-field-boundedness}.
}
\end{proof}

Next, we obtain a quantitative estimate on the continuity of $E$ as below.
\begin{proposition}\label{prop:continuity-of-E}
Let $1/2<\alpha<1$, $0<s_0<s_1<1$, and $\gamma \in C([s_0,s_1];\Leb{2}_{q,p})$. Suppose that 
\begin{equation}\label{eq:E-continuity}
 \partial_s \rho + \Div_q \boldj=0\quad \text{in } [s_0,s_1]\times \mathbb{R}^3,
\end{equation}
where 
\[  \rho(s,q)=\int_{\mathbb{R}^3} \gamma^2(s,q,p)dp\quad \text{and}\quad \boldj(s,q)=\int_{\mathbb{R}^3} p\gamma^2(s,q,p)dp.\]
\begin{enumerate}[label=\textnormal{(\roman*)}]
\item $E=E[\gamma]$ satisfies
	\begin{equation}\label{continuityE}
	\begin{split}
		\norm{E(s_1)-E(s_0)}{L^\infty_q}\apprle  {{(s_1-s_0)^{2\alpha-1}\norm{\bold{j}}{L^{\infty}_{s,q}} +(s_1-s_0)^{2\alpha-1}\norm{\langle p\rangle^2\gamma}{L^{\infty}_{s,q,p}}^2 }}
	\end{split}
	\end{equation}
	\item We also have the corresponding estimate for $\nabla_q E=\nabla_q E[\gamma]:$
	\begin{equation}\label{continuityDE}
		\begin{split}
			&\norm{\nabla_qE(s_1)-\nabla_q E(s_0)}{L^{\infty}_q}\\
			&\apprle   (s_1-s_0)^{2\alpha-1} \norm{\nabla_q\bold{j}}{L^{\infty}_{s,q}}+(s_1-s_0)^{2\alpha-1} [\norm{\langle p\rangle^4\gamma}{L^{\infty}_{s,q,p}}^2 +\norm{\nabla_q\gamma}{L^{\infty}_{s,q,p}}^2 ]
		\end{split}
	\end{equation}
	\item For $\nabla^2_qE[\gamma]$, we also have 
	\begin{equation}\label{continuityDDE}
		\begin{split}
			&\Vert \nabla^2_qE(s_1)-\nabla_q^2E(s_0)\Vert_{L^\infty_q}\\
			&{{\apprle   (s_1 - s_0)^{2\alpha-1}\Vert \nabla_{q}^2{\bf j}\Vert_{L^\infty_{s,q}}
+ (s_1 - s_0)^{2\alpha-1}\left[\Vert 
\langle p\rangle^4\gamma\Vert_{L^\infty_{s,q,p}}\Vert\nabla^2_{q}\gamma\Vert_{L^\infty_{s,q,p}}+\Vert \langle p\rangle^{2}\nabla_q\gamma\Vert_{L^\infty_{s,q,p}}^2\right].}}\\
	\end{split}
	\end{equation}
\end{enumerate} 
\end{proposition}
\begin{proof}
	We only prove (ii) and (iii) since the proofs of  (i) is similar. { Fix $\eta \in C_c^\infty([s_0,s_1))$, $\psi\in C_c^\infty(\mathbb{R}^3)$. Then if we test the equation \eqref{eq:E-continuity} by $v(s,q)=\eta(s)\Phi(q)$, then 
\begin{equation}\label{eq:eta-test-function}
\begin{aligned}
&\int_{s_0}^{s_1}\eta(s)\int_{\mathbb{R}^3} \nabla_q \psi(q) \cdot \boldj(s,q) \myd{q}ds\\
&=-\int_{s_0}^{s_1}\eta'(s)\int_{\mathbb{R}^3} \rho(s,q)\psi(q)\myd{q}ds-\eta(s_0)\int_{\mathbb{R}^3} \psi(q)\rho(s_0,q)\myd{q}.
\end{aligned}
\end{equation}
Since $C_c^\infty([s_0,s_1))$ is dense in $\{\phi \in \Sob{1}{1}([s_0,s_1]) : \phi(s_1)=0\}$, the identity \eqref{eq:eta-test-function} holds for any $\eta \in \Sob{1}{1}([s_0,s_1])$ with $\eta(s_1)=0$. 

Given $s_0\leq t<s_1$ and $0<h<s_1-t$, we take $\eta=\eta_{t,h}$, where
\[
\eta_{t,h}(s)=\begin{dcases}
1 &\quad \text{if } s_0\leq s\leq t,\\
1-\frac{s-t}{h}&\quad \text{if } t\leq s\leq t+h,\\
0&\quad \text{if } t+h\leq s\leq s_1.
\end{dcases}
\]
If we choose $\psi=\eta_{t,h}$ in \eqref{eq:eta-test-function}, then we have \begin{equation}\label{eq:eta-test-function-2}
\begin{aligned}
&\int_{s_0}^{t+h}\eta_{t,h}(s)\int_{\mathbb{R}^3} \nabla_q \psi(q) \cdot \boldj(s,q) \myd{q}ds\\
&=\frac{1}{h}\int_{t}^{t+h}\int_{\mathbb{R}^3} \rho(s,q)\psi(q)\myd{q}ds-\int_{\mathbb{R}^3} \rho(s_0,q)\psi(q)\myd{q}.
\end{aligned}
\end{equation}
Letting $h\rightarrow 0+$, we get 
\begin{equation}\label{eq:eta-test-function-3}
\int_{s_0}^{t}\int_{\mathbb{R}^3} \nabla_q \psi(q) \cdot \boldj(s,q) \myd{q}ds=\int_{\mathbb{R}^3} \rho(t,q)\psi(q)\myd{q}-\int_{\mathbb{R}^3} \rho(s_0,q)\psi(q)\myd{q}.
\end{equation}
Since $\gamma \in C([s_0,s_1];\Leb{2}_{q,p})$, letting $t\rightarrow s_1$, we get 
	}
        \begin{equation}\label{eq:weak-formulation-continuity-eqn}
        \iint_{\mathbb{R}^3_r\times\mathbb{R}^3_u} \psi(r)[\gamma^2(s_1,r,u)-\gamma^2(s_0,r,u)]drdu=\int_{s_0}^{s_1}\iint_{\mathbb{R}^3_r\times\mathbb{R}^3_u} (u\cdot \nabla_r \psi)\gamma^2 \myd{u}drds
                \end{equation}
for all $\psi \in C_c^\infty(\mathbb{R}^3)$. If we put 
\[ \psi(r)=R^{2\alpha-3} \partial_{r^i} \left(\{\partial_{q^j} \chi\}(R^{-1}(q-r))\right), \]
in \eqref{eq:weak-formulation-continuity-eqn}, then the left hand side becomes
\begin{align*}
-\partial_{q^i} E_R^j(s_1)+\partial_{q^i} E_R^j(s_0).
\end{align*}

On the right-hand side, we have 
\begin{align*}
&\int_{s_0}^{s_1} \iint_{\mathbb{R}^3_r\times\mathbb{R}^3_u} (u\cdot \nabla_r [ R^{2\alpha-3} \partial_{r^i} (\{\partial_{q^j} \chi\}(R^{-1}(q-r)))] \gamma^2 dudrds\\
&=\int_{s_0}^{s_1} \int_{\mathbb{R}^3_r} \mathbf{j}\cdot \nabla_r [R^{2\alpha-3} \partial_{r^i} (\{\partial_{q^j} \chi\}(R^{-1}(q-r)))] drds\\
&=-\int_{s_0}^{s_1} \int_{\mathbb{R}^3_r} (\Div \mathbf{j}) [R^{2\alpha-3} \partial_{r^i} (\{\partial_{q^j} \chi\}(R^{-1}(q-r)))]drds\\
&\apprle_\alpha (s_1-s_0) \norm{\nabla_q \mathbf{j}}{\Leb{\infty}_{s,q}} R^{2\alpha-1}.
\end{align*}
This implies that 
\begin{equation}\label{eq:nabla-Lipschitz}
|\nabla E_R(s_1)-\nabla E_R(s_0)|\apprle_\alpha (s_1-s_0) {\norm{\nabla_q\mathbf{j}}{\Leb{\infty}_{s,q}}}R^{2\alpha-1}.
\end{equation}
{By  \eqref{eq:ERV-estimate}, one can show that 
	\begin{equation}\label{eq:nabla-ER}
|\nabla E_R |\apprle_\alpha R^{2\alpha} \norm{\nabla_q\gamma}{L^{\infty}_{s,q,p}}\norm{\langle p\rangle^4\gamma}{L^{\infty}_{s,q,p}}		.
	\end{equation} }
Hence since $\alpha<1$, { it follows from \eqref{eq:nabla-ER}, \eqref{eq:nabla-Lipschitz}, and Young's inequality that }
\begin{align*}
|\nabla E(s_1,q)-\nabla E(s_0,q)|&\leq \int_0^A  |\nabla E_R(s_1,q)-\nabla E_R(s_0,q)|\frac{dR}{R^2}\\
&\relphantom{=}+\int_A^\infty  |\nabla E_R(s_1,q)-\nabla E_R(s_0,q)|\frac{dR}{R^2}\\
&\apprle_\alpha  A^{2\alpha-1}\norm{\nabla_q \gamma}{\Leb{\infty}_{s,q,p}} \norm{\action{p}^4\gamma}{\Leb{\infty}_{s,q,p}}+ (s_1-s_0)\norm{\nabla_q\boldj}{\Leb{\infty}_{s,q}}A^{2\alpha-2}\\
&\apprle_\alpha {A^{2\alpha-1} (\norm{\nabla_q\gamma}{\Leb{\infty}_{s,q,p}}^2+\norm{\action{p}^4\gamma}{\Leb{\infty}_{s,q,p}}^2)+(s_1-s_0)\norm{\nabla_q \boldj}{\Leb{\infty}_{s,q}}A^{2\alpha-2}}.
\end{align*}
for any $A>0$. {Since the implicit constant depends only on $\alpha$, then (ii) follows by choosing $A=(s_1-s_0)$.} 

{ For the proof of (iii), we use the following bound instead of \eqref{eq:nabla-ER}, which can be proved in a similar manner as before, 
\begin{equation*}
	\begin{split}
|\nabla^2 E_R|\lesssim R^{2\alpha}\left[\Vert 
\langle p\rangle^4\gamma\Vert_{L^\infty_{s,q,p}}\Vert\nabla^2_{q}\gamma\Vert_{L^\infty_{s,q,p}}+\Vert \langle p\rangle^{2}\nabla_q\gamma\Vert_{L^\infty_{s,q,p}}^2\right].     		
	\end{split}
\end{equation*}
This completes the proof of Proposition \ref{prop:continuity-of-E}.}
\end{proof}

\begin{proposition}\label{prop:continuity-of-E-2}
Let $1<\alpha<3/2$, $0<s_0<s_1$, and $\gamma \in C([s_0,s_1];\Leb{2}_{q,p})$ satisfy \eqref{eq:E-continuity}. Then we have the following estimate on $E=E[\gamma]$:
\begin{enumerate}[label=\textnormal{(\roman*)}]
	\item $E=E[\gamma]$ satisfies the following estimate:
	\begin{equation}\label{continuityE2}
	\begin{split}
		\norm{E(s_1)-E(s_0)}{L^\infty_q}\lesssim 
		(s_1-s_0)^{2-\alpha}\left(\norm{\bold{j}}{L^{\infty}_{s,q}}+\norm{\gamma}{L^{\infty}_sL^2_{q,p}}^2 \right).
	\end{split}
	\end{equation}
	\item We also have the corresponding estimate for $\nabla_q E=\nabla_q E[\gamma]:$
	\begin{equation}\label{continuityDE2}
		\begin{split}
			\norm{\nabla_qE(s_1)-\nabla_q E(s_0)}{L^{\infty}_q}&\lesssim  (s_1-s_0)^{\frac{5-2\alpha}{3}} \norm{\nabla_q\bold{j}}{L^{\infty}_{s,q}}+(s_1-s_0)^{\frac{5-2\alpha}{3}} \norm{\gamma}{L^{\infty}_s \Leb{2}_{q,p}}^2 .\\
		\end{split}
	\end{equation}
	\item For $\nabla^2_qE[\gamma]$, we also have 
	\begin{equation}\label{continuityDDE2}
		\begin{split}
			\Vert \nabla^2_qE(s_1)-\nabla_q^2E(s_0)\Vert_{L^\infty_q}&\lesssim  (s_1 - s_0)^{\frac{3-\alpha}{2}}\Vert \nabla_{q}^2{\bf j}\Vert_{L^\infty_{s,q}}+ (s_1 - s_0)^{\frac{3-\alpha}{2}}\|\gamma\|_{L^\infty_s L^2_{q,p}}^2.
	\end{split}
	\end{equation}
\end{enumerate}
\end{proposition}
\begin{proof}
{(i) If we put 
\[ \psi(r)=R^{2\alpha-3}(\{\partial_{q^j} \chi\}(R^{-1}(q-r))). \]
in \eqref{eq:weak-formulation-continuity-eqn}, then following a similar argument as in the proof of \eqref{eq:nabla-Lipschitz}, we have
\begin{equation}\label{eq:ER-Lipschitz-alpha-bigger-1}
|E_R(s_1)-E_R(s_0)|\apprle R^{2\alpha-1} \norm{\mathbf{j}}{\Leb{\infty}_{s,q}}.
\end{equation}
Since $\alpha>1$, it follows from \eqref{eq:ER-estimate} and \eqref{eq:ER-Lipschitz-alpha-bigger-1} that for $0<s_0<s_1$, we have 
\begin{align*}
|E(s_1)-E(s_0)|&\leq \int_0^{A^{-1/2}} |E_R(s_1)-E_R(s_0)| \frac{dR}{R^2} + \int_{A^{-1/2}}^\infty |E_R(s_1)-E_R(s_0)| \frac{dR}{R^2}\\
&\apprle_\alpha (s_1-s_0)\norm{\mathbf{j}}{\Leb{\infty}_{s,q}}\int_0^{A^{-1/2}} R^{2\alpha-3} \myd{R}+\norm{\gamma}{\Leb{2}_{s,q,p}}^2\int_{A^{-1/2}}^\infty R^{2\alpha-5} dR\\
&\apprle_\alpha (s_1-s_0) A^{1-\alpha} \norm{\mathbf{j}}{\Leb{\infty}_{s,q}} + A^{2-\alpha} \norm{\gamma}{\Leb{\infty}_s\Leb{2}_{q,p}}^2.
\end{align*}
Hence we get the desired result by choosing $A=(s_1-s_0)$.

(ii) Since $\alpha>1$, it follows from \eqref{eq:nabla-Lipschitz} and \eqref{eq:nabla-ER} that
\begin{align*}
|E(s_1)-E(s_0)|&\leq \int_0^{A^{-1/3}} |E_R(s_1)-E_R(s_0)|\frac{dR}{R^2} +\int_{A^{-1/3}}^\infty |E_R(s_1)-E(s_0)|\frac{dR}{R^2}\\
&\apprle_\alpha (s_1-s_0)\norm{\nabla \boldj}{\Leb{\infty}_{s,q}}\int_0^{A^{-1/3}} R^{2\alpha-3}dR+\norm{\gamma}{\Leb{\infty}_s\Leb{2}_{q,p}}^2\int_{A^{-1/3}}^\infty R^{2\alpha-6}dR\\
&\apprle_\alpha (s_1-s_0)\norm{\nabla\boldj}{\Leb{\infty}_{s,q}}A^{\frac{2-2\alpha}{3}}+\norm{\gamma}{\Leb{\infty}_s\Leb{2}_{q,p}}^2 A^{\frac{5-2\alpha}{3}}.
\end{align*}
Then the desired result follows by choosing $A=(s_1-s_0)$.

(iii) If we put 
\[ \psi(r)=R^{2\alpha-3} \partial_{r^ir^k}(\{\partial_{q^j} \chi\}(R^{-1}(q-r))). \]
in \eqref{eq:weak-formulation-continuity-eqn}, then following a similar argument as in the proof of \eqref{eq:nabla-Lipschitz}, we have
\begin{equation}\label{eq:Hessian-control}
|\nabla_q^2 E_R(s_1)-\nabla_q^2 E_R(s_0)|\apprle_\alpha (s_1-s_0) R^{2\alpha-1}\norm{\nabla_q^2 \mathbf{j}}{\Leb{\infty}_{s,q}}.
\end{equation}
Also, one can show that 
\begin{equation}\label{eq:Hessian-control-2}
|\nabla_q^2 E_R(s)|\apprle R^{2\alpha-5} \norm{\gamma}{\Leb{\infty}_s\Leb{2}_{q,p}}^2,\\
\end{equation}
Since $\alpha>1$, it follows from \eqref{eq:Hessian-control} and \eqref{eq:Hessian-control-2} that 
\begin{align*}
|E(s_1)-E(s_0)|&\leq \int_0^{A^{-1/4}} |E_R(s_1)-E_R(s_0)|\frac{dR}{R^2} +\int_{A^{-1/4}}^\infty |E_R(s_1)-E(s_0)|\frac{dR}{R^2}\\
&\apprle_\alpha (s_1-s_0)\norm{\nabla^2\boldj}{\Leb{\infty}_{s,q}}A^{\frac{2-2\alpha}{4}}+\norm{\gamma}{\Leb{\infty}_s\Leb{2}_{q,p}}^2 A^{\frac{6-2\alpha}{4}}.
\end{align*}
Then the desired result follows by choosing $A=(s_1-s_0)$ in the end. This completes the proof of Proposition \ref{prop:continuity-of-E-2}.
}
\end{proof}

If we assume that $\gamma$ satisfies an equation that is related to the Vlasov-Riesz system, then we get the following estimate on the electric field.

\begin{proposition}\label{prop:quantitative-control-E}
Let $\alpha \in (1/2,3/2)\setminus\{1\}$ and $I\subset [0,1]$ be an interval. If $\gamma$ satisfies
\[
\gamma \in \Leb{\infty}(I;\Leb{\infty}_{q,p})\cap {C(I;\Leb{2}_{q,p})},\quad |p|^2\gamma \in \Leb{\infty}(I;\Leb{\infty}_{q,p}),
\] 
and
\begin{equation}\label{eq:transport}
\partial_s \{ \gamma^2 \} + \Div_q \{ p \gamma^2\} +\Div_p \{F \gamma^2 \}=0\quad \text{in } I\times \mathbb{R}^3_q\times \mathbb{R}_p^3
\end{equation}
for some force field $F(s,q)$, then 
\begin{enumerate}[label=\textnormal{(\roman*)}]
	\item If $1/2<\alpha<1$, then
\begin{align*}
\norm{E[\gamma](s_1)-E[\gamma](s_0)}{\Leb{\infty}_q}&\apprle (s_1-s_0)^{2\alpha-1}\action{\ln(s_1-s_0)}\norm{|p|^2\gamma}{\Leb{\infty}_s(I;\Leb{\infty}_{q,p})}^2\\
&\relphantom{=}+ (s_1-s_0)^{\beta_2}\norm{\action{p}^2\gamma}{\Leb{\infty}_s(I;\Leb{\infty}_{q,p})}^2\\
&\relphantom{=}+(s_1-s_0)^{2\alpha-1} \left[\norm{\gamma}{\Leb{\infty}_s(I;\Leb{\infty}_{q,p})}^2+\norm{\gamma}{\Leb{\infty}_s(I;\Leb{2}_{q,p})}^2\right]\\
&\relphantom{=}+(s_1-s_0)^{\beta_3}\norm{\action{p}^2\gamma}{\Leb{\infty}_s(I;\Leb{\infty}_{q,p})}^2\norm{s^{2\alpha-1}F}{\Leb{\infty}_{s}(I;\Leb{\infty}_{q})}
\end{align*}

\item If $1<\alpha<3/2$, then
\begin{align*}
\norm{E[\gamma](s_1)-E[\gamma](s_0)}{\Leb{\infty}_q}&\apprle (s_1-s_0)^{2-\alpha}\action{\ln(s_1-s_0)}\norm{|p|^2\gamma}{\Leb{\infty}_s(I;\Leb{\infty}_{q,p})}^2\\
&\relphantom{=}+ (s_1-s_0)^{\beta_2}\norm{\action{p}^2\gamma}{\Leb{\infty}_s(I;\Leb{\infty}_{q,p})}^2\\
&\relphantom{=}+(s_1-s_0)^{2-\alpha} \left[\norm{\gamma}{\Leb{\infty}_s(I;\Leb{\infty}_{q,p})}^2+\norm{\gamma}{\Leb{\infty}_s(I;\Leb{2}_{q,p})}^2\right]\\
&\relphantom{=}+s_0^{2-2\alpha}(s_1-s_0)^{\beta_3}\norm{\action{p}^2\gamma}{\Leb{\infty}_s(I;\Leb{\infty}_{q,p})}^2\norm{s^{2\alpha-1}F}{\Leb{\infty}_{s}(I;\Leb{\infty}_{q})}
\end{align*}
\end{enumerate}
for all $s_0,s_1 \in I$ satisfying $0<s_0<s_1$, where 
{\begin{equation}\label{eq:beta-order}
\beta_2,\beta_3 > \max\{10\alpha-10,2\alpha-1\}+10-9\alpha
\end{equation}}
\end{proposition}

\begin{proof}
{ Since $\gamma$ is a solution of \eqref{eq:transport}, following a cut-off argument as in the proof of Proposition \ref{prop:continuity-of-E}, one can show that $\gamma$ satisfies 
\begin{equation}\label{eq:cut-off}
\begin{aligned}
0&=\iint_{\mathbb{R}^3_q\times\mathbb{R}^3_p} \gamma^2(s_1,q,p)\myd{q}dp-\iint_{\mathbb{R}^3_q\times\mathbb{R}^3_p} \gamma^2(s_0,q,p) \myd{q}dp\\
&\relphantom{=}+\int_{s_0}^{s_1}\iint_{\mathbb{R}^3_q\times\mathbb{R}^3_p}  p\gamma^2 \cdot \nabla_q\phi \myd{q}dpds+\int_{s_0}^{s_1}\iint_{\mathbb{R}^3_q\times\mathbb{R}^3_p}  (F\cdot \nabla_p \phi)\gamma^2  \myd{q}dpds
\end{aligned}
\end{equation}
for all $\phi \in C_c^\infty(\mathbb{R}^3_q\times \mathbb{R}^3_p)$ and $s_0,s_1 \in I$ with $s_0<s_1$.}

(i) For simplicity, we write $E=E[\gamma]$. {{Fix $s,s_1,s_0\in I$. Let $0<A,B<1$ to be chosen later.}}   Define $\theta = (2\alpha-1)/(2\alpha-4)<0$. {{Similar to \eqref{eq:E-representation-V-R}, we have
\begin{equation*}
	\begin{aligned}
		E(s,q)&=\int_0^A E_R\frac{dR}{R^2}+\int_{A^\theta}^\infty E_R\frac{dR}{R^2}+\int_0^B\int_A^{A^\theta}E_{R,V}\frac{dR}{R^2}\frac{dV}{V}+\int_B^{B^{-3}}\int_A^{A^\theta}E_{R,V}\frac{dR}{R^2}\frac{dV}{V}\\
		&\quad +\int_{B^{-3}}^\infty\int_A^{A^\theta}E_{R,V}\frac{dR}{R^2}\frac{dV}{V}.
	\end{aligned}
\end{equation*}
}}{{Then using \eqref{eq:ER-estimate} and \eqref{eq:ER-estimate-3}}}, it is easy to see that 
\begin{equation}\label{eq:ER-estimate-2}
\begin{split}
\int_{A^{\theta}}^\infty |E_R(s)|\frac{dR}{R^2}&\apprle_\alpha A^{2\alpha-1} \norm{\gamma}{\Leb{2}_{q,p}}^2,\\
\int_0^A |E_R(s)|\frac{dR}{R^2}&\apprle_\alpha A^{2\alpha-1} \norm{|p|^2 \gamma}{\Leb{\infty}_{q,p}}^2.
\end{split}
\end{equation}
Also, observe that  \eqref{eq:ERV-estimate} gives
\begin{equation}\label{eq:ERV-estimate-2}
\begin{aligned}
\int_0^{A^\theta} \int_0^B |E_{R,V}(s)| \frac{dV}{V}\frac{dR}{R^2} &\apprle_\alpha A^{\theta(2\alpha-1)} B^3 \norm{\gamma(s)}{\Leb{\infty}_{q,p}}^2,\\
\int_0^{A^\theta} \int_{B^{-3}}^\infty |E_{R,V}(s)| \frac{dV}{V}\frac{dR}{R^2} &\apprle_\alpha A^{\theta(2\alpha-1)} B^3 \norm{\action{p}^2\gamma(s)}{\Leb{\infty}_{q,p}}^2.
\end{aligned}
\end{equation}
It remains to estimate
\begin{equation*}
	\begin{aligned}
	\int_A^{A^\theta}\int_B^{B^{-3}}(E_{R,V}(s_1)-E_{R,V}(s_2))\frac{dV}{V}\frac{dR}{R^2}.	
	\end{aligned}
\end{equation*}

Note that 
\begin{equation*}
	\begin{aligned}
	&\int_A^{A^\theta}\int_B^{B^{-3}}	E_{R,V}(s,q)\frac{dV}{V}\frac{dR}{R^2}\\
	&=\int_A^{A^\theta} \left[\iint_{\mathbb{R}^3_{r}\times\mathbb{R}^3_{u}}  R^{-(3-2\alpha)} (\nabla \chi)(R^{-1}(q-y))\left(\int_B^{B^{-3}} \chi(V^{-1}u)\frac{dV}{V}\right)\gamma^2(s,y,u)dudy\right]\frac{dR}{R^2}\\
		&:=\int_A^{A^\theta} \mathcal{E}_{R,a}(s,q)\frac{dR}{R^2}.
	\end{aligned}
\end{equation*}
where 
\[  \mathcal{E}_{R,a}(s,q):=\iint_{\mathbb{R}^3_{r}\times\mathbb{R}^3_{u}}  R^{-(3-2\alpha)} \{\nabla \chi\} (R^{-1}(q-r))\cdot \chi_{\{B\leq \cdot \leq B^{-3}\}}(u)\cdot \gamma^2(r,u)drdu.\]
and
\[  \chi_{\{B\leq \cdot \leq B^{-3}\}} (u) =\int_{\{B\leq V\leq B^{-3}\}} \chi(V^{-1} u)\frac{dV}{V}.\]
On the other hand, if we choose $\phi(r,u)=R^{-(3-2\alpha)}\partial_{q^j} \chi(R^{-1}(q-r))\cdot \chi_{\{B\leq \cdot \leq B^{-3}\}}(u)$ in \eqref{eq:cut-off}, then we have
\begin{align*}
0&=\mathcal{E}_{R,a}^j(s_1,q)-\mathcal{E}_{R,a}^j(s_0,q)\\
&\relphantom{=}-\int_{s_0}^{s_1} \iint_{\mathbb{R}^3_{r}\times\mathbb{R}^3_{u}}  R^{-(4-2\alpha)} u^k \{\partial_{q^j}\partial_{q^k} \chi\} (R^{-1}(q-r))\cdot \chi_{\{B\leq \cdot \leq B^{-3}\}}(u)\cdot \gamma^2(s,r,u)drduds \\
&\relphantom{=}+\int_{s_0}^{s_1} \iint_{\mathbb{R}^3_{r}\times\mathbb{R}^3_{u}}  R^{-(3-2\alpha)} \partial_{q^j} \chi(R^{-1}(q-r)) \gamma^2  F\cdot (\nabla_u \chi_{\{B\leq \cdot \leq B^{-3}\}} )(u)drduds.
\end{align*}
Since $\supp \chi \subset B_2 \setminus B_{1/2}$, one can show that
\[  |\nabla_u \chi_{\{B\leq \cdot \leq B^{-3}\}}(u)|\apprle B^{-1} 1_{\{|u|\leq 2B\}} + B^3 1_{\{|u|\geq B^{-3}/2\}}.\]
Then we have
\begin{align*}
&\left|\iint_{\mathbb{R}^3_{r}\times\mathbb{R}^3_{u}} R^{-(3-2\alpha)} \partial_{q^j} \chi(R^{-1}(q-r))\cdot \gamma^2 (r,u)\cdot (F\cdot\nabla_u)\chi_{\{B\leq \cdot \leq B^{-3}\}}(u)drdu \right|\\
&\apprle_\alpha \norm{s^{2\alpha-1}F}{\Leb{\infty}_q} \cdot s^{1-2\alpha} R^{2\alpha} \cdot \left[ B^2 \norm{\gamma}{\Leb{\infty}_{q,p}}^2 + B^6 \norm{|p|^2 \gamma}{\Leb{\infty}_{q,p}}^2\right].
\end{align*}
Hence it follows that 
\begin{align*}
&\left|\int_{A}^{A^{\theta}} \{\mathcal{E}_{R,a}(s_1,q)-\mathcal{E}_{R,a}(s_0,q)\}\frac{dR}{R^2}\right|\\
&\apprle_\alpha \left(\int_{s_0}^{s_1} s^{1-2\alpha} ds\right)\left(\int_A^{A^\theta} \frac{dR}{R^{2-2\alpha}}\right) B^2\norm{\action{p}^2\gamma}{\Leb{\infty}(I;\Leb{\infty}_{q,p})}^2\norm{s^{2\alpha-1}F}{\Leb{\infty}(I;\Leb{\infty}_q)}\\
&\relphantom{=}+\int_A^{A^\theta} \int_{s_0}^{s_1} \iint_{\mathbb{R}^3_{r}\times\mathbb{R}^3_{u}} R^{-(4-2\alpha)} |u||\nabla^2 \chi|(R^{-1}(q-r)) \left(\int_{B}^{B^{-3}} \frac{\chi(V^{-1} u)}{V} dV \right)\gamma^2(s,r,u) drduds \frac{dR}{R^2}\\
&\apprle_\alpha \left(\int_{s_0}^{s_1} s^{1-2\alpha} ds\right)\left(\int_A^{A^\theta} \frac{dR}{R^{2-2\alpha}}\right)  B^2\norm{\action{p}^2\gamma}{\Leb{\infty}(I;\Leb{\infty}_{q,p})}^2 \norm{s^{2\alpha-1}F}{\Leb{\infty}(I;\Leb{\infty}_q)}\\
&\relphantom{=}+(s_1-s_0)\norm{|p|^2\gamma}{\Leb{\infty}(I;\Leb{\infty}_{q,p})}^2\int_A^{A^\theta} \frac{dR}{R^{3-2\alpha}} \int_{B}^{B^{-3}} \frac{dV}{V}\\
&\apprle_\alpha |s_1^{2-2\alpha}-s_0^{2-2\alpha}|A^{\theta(2\alpha-1)}B^2\norm{\action{p}^2\gamma}{\Leb{\infty}(I;\Leb{\infty}_{q,p})}^2 \norm{s^{2\alpha-1}F}{\Leb{\infty}(I;\Leb{\infty}_q)}\\
&\relphantom{=}+(s_1-s_0)\norm{|p|^2\gamma}{\Leb{\infty}(I;\Leb{\infty}_{q,p})}^2{{A^{\theta(2\alpha-2)}}}\ln (B^{-4}).
\end{align*}
Combining this with \eqref{eq:ER-estimate-2} and \eqref{eq:ERV-estimate-2}, we get 
\begin{align*}
|E(s_1)-E(s_0)|&\apprle_\alpha A^{\theta (2\alpha-1)}B^3\norm{\action{p}^2\gamma}{\Leb{\infty}(I;\Leb{\infty}_{q,p})}^2 + A^{2\alpha-1}(\norm{\action{p}^2\gamma}{\Leb{\infty}(I;\Leb{\infty}_{q,p})}^2+\norm{\gamma}{\Leb{\infty}(I;\Leb{2}_{q,p})}^2)
 \\
&\relphantom{=}+ |s_1^{2-2\alpha}-s_0^{2-2\alpha}| A^{\theta(2\alpha-1)} B^2\norm{\action{p}^2\gamma}{\Leb{\infty}(I;\Leb{\infty}_{q,p})}^2\norm{s^{2\alpha-1}F}{\Leb{\infty}(I;\Leb{\infty}_q)}\\
&\relphantom{=}+(s_1-s_0)A^{\theta(2\alpha-2)}\ln (B^{-4}) \norm{|p|^2 \gamma}{\Leb{\infty}(I;\Leb{\infty}_{q,p})}^2.
\end{align*}
Finally, if we choose
\[ A=(s_1-s_0)\quad \text{and}\quad B=(s_1-s_0)^\beta,\]
then
\begin{align*}
|E(s_1)-E(s_0)|&\apprle_\alpha (s_1-s_0)^{\theta(2\alpha-1)+3\beta}\norm{\action{p}^2\gamma}{\Leb{\infty}(I;\Leb{\infty}_{q,p})}^2\\
&\relphantom{=} +(s_1-s_0)^{2\alpha-1} \left(\norm{\action{p}^2\gamma}{\Leb{\infty}(I;\Leb{\infty}_{q,p})}^2+\norm{\action{p}^2\gamma}{\Leb{\infty}(I;\Leb{2}_{q,p})}^2\right)\\
&\relphantom{=}+ |s_1^{2-2\alpha}-s_0^{2-2\alpha}| (s_1-s_0)^{\theta(2\alpha-1)+2\beta}  \norm{\action{p}^2 \gamma}{\Leb{\infty}(I;\Leb{\infty}_{q,p})}^2\norm{s^{2\alpha-1}F}{\Leb{\infty}(I;\Leb{\infty}_q)}\\
&\relphantom{=}+(s_1-s_0)^{2\alpha-1}\action{\ln(s_1-s_0)} \norm{|p|^2 \gamma}{\Leb{\infty}(I;\Leb{\infty}_{q,p})}^2.
\end{align*}
where we used $1+\theta(2\alpha-2)>2\alpha-1$ in the last line. 
Hence, the desired result follows by choosing $\beta$ sufficiently large. 

(ii) The proof is similar to that of (i). We only highlight the key difference. Similarly using \eqref{eq:ER-estimate} and \eqref{eq:ERV-estimate}, if we choose $\theta=-1/2, \beta=\frac{2-\alpha}{2\alpha-1}$, then 
\[
\int_{A^{\theta}}^\infty |E_R(s)|\frac{dR}{R^2}\apprle_\alpha A^{2-\alpha} \norm{\gamma}{\Leb{2}_{q,p}}^2\quad \text{and}\quad 
\int_0^{A^\beta} |E_R(s)|\frac{dR}{R^2}\apprle_\alpha A^{2-\alpha} \norm{\action{p}^2 \gamma}{\Leb{\infty}_{q,p}}^2.
\]
Then the proof is exactly the same by replacing $[A,A^\theta]$ with $[A^\beta,A^\theta]$ in each step. We also note that 
\[ {{|s_1^{2-2\alpha}}}-s_0^{2-2\alpha}|\leq 2s_0^{2-2\alpha} \]
if $\alpha>1$. This completes the proof of Proposition \ref{prop:quantitative-control-E}.
\end{proof}

\section{Proof of Theorem \ref{thm:A}}\label{sec:3}
This section is devoted to proving Theorem \ref{thm:A}. The following lemma is an immediate application of \eqref{eq:estimate-transport} and \eqref{eq:commutator-relationship}. This lemma is necessary for us to perform bootstrap argument. 

\begin{lemma}\label{lem:transport}
Let  $\alpha\in (1/2,3/2)\setminus\{1\}$ and $\gamma$ be a strong solution of \eqref{eq:VR-pct} on $T^*\leq s\leq 1$ with ``initial'' data $\gamma(s=1)=\gamma_1$. Assume that $\gamma_1$ satisfies for some $a\in \mathbb{N}$, $r\in [2,\infty]$ that 
\[ \norm{\action{p}^a \gamma_1}{\Leb{r}_{q,p}}\leq \varepsilon_0\]
and that 
\[   |E(s,q)|\leq D,\quad T^*\leq s\leq 1.\]
Then there holds that 
\[ \norm{\gamma(s)}{\Leb{r}_{q,p}}\leq \varepsilon_0\quad\text{and}\quad {{\norm{|p|^a \gamma(s)}{\Leb{r}_{q,p}}\apprle }}\varepsilon_0(1+aD)s^{a(2-2\alpha)}.\]
In particular, 
\[ \norm{\action{p}^a\gamma(s)}{L^r_{q,p}}\lesssim\varepsilon_0(1+aD)\max{\{1, s^{a(2-2\alpha)} \}}. \]
\end{lemma} 
{ \begin{proof}
Since $\mathcal{L}[\gamma]=0$, using \eqref{eq:estimate-transport},	one can immediately get the first estimate. For the second estimates, we first consider the case when $a=1$. Using \eqref{eq:commutator-relationship}, we have $\mathcal{L}_\gamma[p^j\gamma]=\lambda s^{1-2\alpha} E^j[\gamma]\gamma$. Using \eqref{eq:estimate-transport}, we can get for any $j=1,2,3$,
\begin{equation*}
	\begin{aligned}
	\norm{p^j\gamma(s)}{L^r_{q,p}}\leq \varepsilon_0+\int_s^1\varepsilon_0\lambda \tau^{1-2\alpha}\norm{E^j}{L^\infty_{s,q}}d\tau\leq \varepsilon_0(1+aD)s^{2-2\alpha}	.
	\end{aligned}
\end{equation*} 
By induction on $a$, we can get the desired estimates.
\end{proof}}

{
\begin{remark}
In the case when $1/2<\alpha<1$, we have $0<s^{2-2\alpha}<1$. So, the moments $\action{p}^a\gamma$ are bounded in the Lebesgue space. In the case when $\alpha>1$, we have $s^{2-2\alpha}>1$. Hence the norm of moments would grow polynomially as $s$ goes to 0. This is the essential reason why we have modified scattering when $\alpha>1$ but scattering when $1/2<\alpha<1$.
\end{remark}
}

The following proposition is a key tool to prove the global well-posedness and scattering of solutions.

\begin{proposition}\label{prop:bootstrap}
There exists an absolute constant $\delta>0$ such that the following holds when $\alpha \in (1/2,1+\delta)\setminus\{1\}$: { there exists a small $\varepsilon_1>0$ such that if $\gamma$ is a solution of \eqref{eq:VR-pct} on $T^*\leq s\leq 1$ with ``initial'' data $\gamma(s=1)=\gamma_1$ satisfying 
\begin{equation}\label{eq:initial-control}
 \norm{\gamma_1}{\Leb{2}_{q,p}}+\norm{\gamma_1}{\Leb{\infty}_{q,p}}\leq \varepsilon_1,
 \end{equation}
 then we have the following bootstrap results:
 }

\begin{enumerate}[label=\textnormal{(\roman*)}]
\item (Moments and the electric field) if there holds that 
\begin{equation}\label{eq:initial-control-2}
 \norm{\action{p}\gamma_1}{\Leb{2}_{q,p}}+\norm{\action{p}^m \gamma_1}{\Leb{\infty}_{q,p}}\leq \varepsilon_1,\quad m\geq 2
\end{equation}
then the electric field $E(s)$ remains bounded and the solution satisfies the bounds 
{{\begin{equation}\label{eq:p-gamma-estimate-l2}
 \norm{\action{p}\gamma(s)}{\Leb{2}_{q,p}}\apprle \varepsilon_1 \action{s^{2-2\alpha}}
 \end{equation}
and
\begin{equation}\label{eq:p-gamma-estimate-l-infty}
\norm{\action{p}^a\gamma(s)}{\Leb{\infty}_{q,p}}\apprle \varepsilon_1 \action{s^{a(2-2\alpha)}},\quad 0\leq a\leq m.
\end{equation}
}}Moreover, there exist $C>0$, $\beta^\prime>0$ such that
\[
  |E(s_1,q)-E(s_2,q)|\leq C\varepsilon_1^2 (s_2-s_1)^{\beta^\prime}\action{\ln(s_2-s_1)}
\]
for all $q\in \mathbb{R}^3$ and $T^*\leq s_1<s_2\leq 1$. {{Moreover, $\beta^\prime=2\alpha-1$ for $1/2<\alpha<1$ and $\beta^\prime=10-9\alpha$ for $1<\alpha<3/2$.}}
\item (Derivatives) Assume additionally that for some $b\in\{0,1\}$, there holds that 
$$ \|\langle p\rangle^b \nabla_{p,q}\gamma_1\|_{L^{\infty}_{q,p}}\leq \varepsilon_1. 
$$
Then we have
\begin{equation}\label{eq:weight-p-a}
\begin{aligned}
	\|\langle p\rangle^a\nabla_{p}\gamma(s)\|_{L^{\infty}_{q,p}}&\apprle \varepsilon_1 \action{s^{a(2-2\alpha)}}  \qquad 0\leq a\leq b \\
	\|\langle p\rangle^a\nabla_{q}\gamma(s)\|_{L^{\infty}_{q,p}}&\apprle \varepsilon_1 h_{\alpha,a}(s)\qquad 0\leq a\leq b
\end{aligned}
\end{equation}
where 
\[  h_{\alpha,a}(s)=\begin{cases}
1 & \quad \text{if } \frac{1}{2}<\alpha<1,a=0,1,\\
s^{\theta_1}&\quad \text{if } \alpha>1, a=0,\\
s^{\theta_2}&\quad \text{if } \alpha>1, a=1
\end{cases}
\]
for some $1-2\alpha<\theta_2<\theta_1<0$ and $\theta_1>-1$. 
\end{enumerate} 
\end{proposition}
\begin{proof}
{(i) Let $C>2$ be a constant larger than all the implied constants appear in the previous section and let $\varepsilon_1>0$ be small enough so that 
\[
4C^2_* \varepsilon_1^2\leq 1,
\]
where $C_*\geq C$ and $\varepsilon_1$ will be determined later.

Define 
\[  I = \{ s\in [T^*,1] : \norm{E(s)}{\Leb{\infty}_q}\leq 2C^2_* \varepsilon_1^2\}.\]
By \eqref{eq:E-estimate} and \eqref{eq:initial-control-2}, we have 
\[ \norm{E(1)}{\Leb{\infty}_q}\leq C(\norm{\gamma_1}{\Leb{2}_{q,p}}^2+\norm{\action{p}^2\gamma_1}{\Leb{\infty}_{q,p}}^2)\leq  2C\varepsilon_1^2. \]
So $1\in I\neq \varnothing$. By continuity, $I$ is closed. It remains for us to show that $I$ is open in $[T^*,1]$. By Lemma \ref{lem:transport}, for $s\in I$, we have 
\begin{equation}\label{eq:control-momentum-gamma}
\begin{aligned}
\norm{\action{p}^a\gamma(s)}{\Leb{r}_{q,p}}&\leq C_0\varepsilon_1(1+a\norm{E}{\Leb{\infty}_{s,q}})s^{a(2-2\alpha)}\\
&\leq C_0\varepsilon_1(1+2aC^2_*\varepsilon_1^2)s^{a(2-2\alpha)}
\end{aligned}
\end{equation}
for $0\leq a\leq m$ and $r\in \{2,\infty\}$.

If $1/2<\alpha<1$, then for $s_1,s_2\in I$ satisfying $0<s_1<s_2\leq 1$, it follows from Proposition \ref{prop:quantitative-control-E} and \eqref{eq:control-momentum-gamma} that 
\begin{equation}\label{eq:convergence-rate}
\begin{aligned}
&\norm{E(s_1)-E(s_2)}{\Leb{\infty}_q}\\
&\leq C (s_2-s_1)^{2\alpha-1}\action{\ln(s_2-s_1)} \norm{|p|^2\gamma}{\Leb{\infty}(I;\Leb{\infty}_{q,p})}^2 \\
&\relphantom{=}+C(s_2-s_1)^{\beta_2} \norm{|p|^2\gamma}{\Leb{\infty}(I;\Leb{\infty}_{q,p})}^2 \\
&\relphantom{=}+C(s_2-s_1)^{2\alpha-1} \left[\norm{\gamma}{\Leb{\infty}(I;\Leb{\infty}_{q,p})}^2 +  \norm{\gamma}{\Leb{\infty}(I;\Leb{2}_{q,p})}^2 \right]\\
&\relphantom{=}+C(s_2-s_1)^{\beta_3} {\norm{\action{p}^2\gamma}{\Leb{\infty}(I;\Leb{\infty}_{q,p})}^2} \norm{E}{\Leb{\infty}(I;\Leb{\infty}_q)}\\
&\leq C(s_2-s_1)^{2\alpha-1}\action{\ln(s_2-s_1)}(C_0^2\varepsilon_1^2(1+4C^2_*\varepsilon_1^2)^2(2+2C^2_*\varepsilon_1^2)+2\varepsilon_1^2)\\
&\leq C(10C_0^2+2)\varepsilon_1^2 (s_2-s_1)^{2\alpha-1}\action{\ln(s_2-s_1)}.
\end{aligned}
\end{equation}
In particular, we have 
\begin{equation}\label{eq:electric-field-difference}
 \norm{E(s_1)-E(s_2)}{\Leb{\infty}_q}\leq C_1 \varepsilon_1^2\frac{\action{k}^2}{2^{k(2\alpha-1)}}
\end{equation}
for $2^{-k}\leq s_1\leq s_2\leq 2^{1-k}$, $k\geq 1$,  and $s_1,s_2\in I$, 
where we chose $C_1>0$ so that
\[ C(10C_0^2+2)\leq C_1.\]
Since $s\in I$ and $1/2<\alpha<1$, it follows from \eqref{eq:electric-field-difference} that 
\begin{equation}
\norm{E(s)}{\Leb{\infty}_q} \leq \norm{E(1)}{\Leb{\infty}_q}+C_1\varepsilon_1^2 \sum_{k=1}^\infty \frac{\action{k}^2}{2^{k(2\alpha-1)}}\leq (2C+C_2)\varepsilon_1^2.
\end{equation}

Now we first choose $C_*\geq 2C+C_2$ and then choose $\varepsilon_1>0$ so that $4C_*^2\varepsilon_1^2\leq 1$. This implies that  $I$ is open in $[T^*,1]$.

If $1<\alpha<10/9$, then for $s_1,s_2\in I$ satisfying $0<s_1<s_2\leq 1$, then it follows from Proposition \ref{prop:quantitative-control-E} and \eqref{eq:control-momentum-gamma} that 
\begin{align*}
\norm{E(s_1)-E(s_2)}{\Leb{\infty}_q}&\leq C(s_2-s_1)^{2-\alpha}\action{\ln(s_2-s_1)}s^{4(2-2\alpha)}_1C_0^2\varepsilon_1^2(1+4C_*^2\varepsilon_1^2)^2\\
&\relphantom{=}+C(s_2-s_1)^{\beta_2}s_1^{4(2-2\alpha)}C_0^2\varepsilon_1^2(1+4C_*^2\varepsilon_1^2)^2\\
&\relphantom{=}+2C\varepsilon_1^2(s_2-s_1)^{2\alpha-1}\\
&\relphantom{=}+C\varepsilon_1^2 s_1^{5(2-2\alpha)}(s_2-s_1)^{\beta_3}(C_0(1+4C_*^2\varepsilon_1^2))^2(1+4C_*^2\varepsilon_1^2)^2\\
&\leq 4C_0^2 C\varepsilon_1^2(s_2-s_1)^{2-\alpha}\action{\ln(s_2-s_1)}s_1^{4(2-2\alpha)}  \\
&\relphantom{=}+4 C_0^2C\varepsilon_1^2(s_2-s_1)^{\beta_2} s_1^{4(2-2\alpha)}+2C\varepsilon_1^2(s_2-s_1)^{2\alpha-1}\\
&\relphantom{=}+16C_0^2 C \varepsilon_1^2 s_1^{5(2-2\alpha)}(s_2-s_1)^{\beta_3}.
\end{align*}
In particular, we have
\begin{equation}\label{eq:E1-alpha-1}
\begin{aligned}
&\norm{E(s_1)-E(s_2)}{\Leb{\infty}_q}\\
&\leq C_3 \varepsilon_1^2\left[2^{-(2-\alpha+8-8\alpha)k}(1+k^2)+2^{-(\beta_2+8-8\alpha)k}+2^{-(2\alpha-1)k} +2^{-k(10-10\alpha+\beta_3)} \right]
\end{aligned}
\end{equation}
for $2^{-k}\leq s_1\leq s_2\leq 2^{1-k}$. Since $1<\alpha<10/9$ and \eqref{eq:beta-order} holds, we have
\begin{equation}\label{eq:growing-rate}
2-\alpha+8-8\alpha=10-9\alpha>0,\quad \beta_2+8-8\alpha>10-9\alpha,\quad 10-10\alpha+\beta_3>10-9\alpha.
\end{equation}

Since $s\in I$, it follows from  \eqref{eq:E1-alpha-1} that 
\begin{align*}
& \norm{E(s)}{\Leb{\infty}_q}\\
 &\leq \norm{E(1)}{\Leb{\infty}_q}+C_3\varepsilon_1^2 \sum_{k=1}^\infty\left[2^{-(2-\alpha+8-8-\alpha)k}(1+k^2)+2^{-(\beta_2+8-8\alpha)k}+2^{-(2\alpha-1)k} +2^{-k(10-10\alpha+\beta_3)} \right]\\
&\leq (2C+C_4)\varepsilon_1^2.
\end{align*}
Now we first choose $C_*\geq 2C+C_4$ and then choose $\varepsilon_1>0$ so that $4C_*^2\varepsilon_1^2\leq 1$. This implies that $I$ is open in $[T^*,1]$, which completes the proof of (i). }

To prove (ii), we use a similar bootstrap argument. We only prove the case $\alpha>1$ since the proof of the remaining case $1/2<\alpha<1$ is similar and, in fact, much simpler since we do not have a growth of $\gamma$ in the weighted $\Leb{\infty}$ spaces in time variable $s$. Suppose that 
\begin{equation}\label{eq:bootstrap-p-moment}
\left\|\langle p\rangle^b \nabla_p \gamma(s)\right\|_{L_{q, p}^{\infty}}  \leq 2 C^4 \varepsilon_1 s^{b(2-2\alpha)}\quad\text{and}\quad \left\|\langle p\rangle^b \nabla_q \gamma(s)\right\|_{L_{q, p}^{\infty}} \leq 2 C^2 \varepsilon_1 h_{\alpha,b}(s).
\end{equation}
In the following, for simplicity, we choose $\varepsilon_1C^4\ll 1$ so that we often suppress the big constant $C$ by $\varepsilon_1$.
Then {{by  \eqref{eq:estimate-transport}, \eqref{eq:commutator-relationship}, \eqref{eq:initial-control-2} and \eqref{eq:bootstrap-p-moment}}} , we have
\begin{align*}
\|\nabla_p \gamma(s) \|_{L^{\infty}_{q,p}}&\leq\|\nabla_p\gamma(1) \|_{L^{\infty}_{q,p}}+\int_s^1 \|\nabla_q\gamma(\tau) \|_{L^{\infty}_{q,p}} d\tau\leq \varepsilon_1+ 2C^2 \varepsilon_1\int_s^1 \tau^{\theta_1} d\tau\leq 2C^3\varepsilon_1 .
\end{align*}
{{Next, suppose that $1<\alpha<1+\delta$, where $\delta>0$ will be determined later. Then by \eqref{eq:estimate-transport}, \eqref{eq:commutator-relationship}, \eqref{eq:gradient-control-less-than-2}  and \eqref{eq:bootstrap-p-moment}, we have
\begin{align*}
\norm{\nabla_q \gamma(s)}{\Leb{\infty}_{q,p}}&\apprle \norm{\nabla_q \gamma(1)}{\Leb{\infty}_{q,p}}+\varepsilon_1\int_s^1 \norm{\nabla_q E}{\Leb{\infty}_q}\tau^{1-2\alpha}d\tau\\ 
&\leq C_0 \varepsilon_1 +\varepsilon_1^3 s^{\beta(5-2\alpha)+2-2\alpha} + \varepsilon_1^3 s^{\beta(2-2\alpha-\theta)+10-10\alpha}+\varepsilon_1^2 \int_s^1 \tau^{\beta(1-2\alpha+3\theta) +1-2\alpha} \norm{\nabla_q\gamma}{\Leb{\infty}_q} d\tau.
\end{align*}
Since $1-2\alpha+3\theta<0$, if we choose $\beta<0$ so that $\beta(1-2\alpha+3\theta)+1-2\alpha>-1$, then it follows from Gronwall's inequality that
\begin{align*}
\norm{\nabla_q \gamma}{\Leb{\infty}_{q,p}}&\leq C_0\varepsilon_1 (1+s^{\beta(5-2\alpha)+2-2\alpha} + s^{\beta(2-2\alpha-\theta)+10-10\alpha}).
\end{align*}
Since $2\alpha-2+\theta >0$, for $\alpha$ sufficiently close to $1$, we have 
\[ -1<\theta_1=\min\{\beta(5-2\alpha)+2-2\alpha,\beta(2-2\alpha-\theta)+10-10\alpha\}<0.\]
Hence we close the bootstrap since we choose a sufficiently large $C$ so that $C\gg C_5$:
\begin{equation}\label{eq:nabla-q-gamma-2}
\norm{\nabla_q \gamma}{\Leb{\infty}_{q,p}}\leq C_5 \varepsilon_1 s^{\theta_1} .
\end{equation}
}}
A similar argument shows that 
\begin{align*}
	\left\|p^m \partial_{p^n} \gamma(s)\right\|_{L_{q, p}^{\infty}} &\leq  \left\||p| \nabla_p \gamma(1)\right\|_{L_{q, p}^{\infty}} +\int_s^1\left\||p| \nabla_q \gamma\left(\tau\right)\right\|_{L_{q, p}^{\infty}} \myd \tau+\int_s^1\|E\|_{L_q^{\infty}}\left\|\nabla_p \gamma\left(\tau\right)\right\|_{L_{q, p}^{\infty}} \tau^{1-2\alpha} \myd\tau\\
	& \leq \varepsilon_1+ \varepsilon_1\int_s^1 \tau^{\theta_2}d\tau + \varepsilon_1^2\int_s^1 \tau^{1-2\alpha}d\tau\apprle \varepsilon_1 s^{2-2\alpha}
\end{align*}

It remains for us to show that 
\[ \norm{\action{p} \nabla_q \gamma(s)}{\Leb{\infty}_{q,p}}\apprle \varepsilon_1 s^{\theta_2}\]
for some $\theta_2<0$.
By \eqref{eq:gradient-control-less-than-2} and \eqref{eq:nabla-q-gamma-2}, we have {
\begin{equation}\label{eq:nabla-q-E-epsilon-2}
\begin{aligned}
 \norm{\nabla_q E}{\Leb{\infty}}&\apprle s^{\beta(5-2\alpha)}\norm{\gamma}{\Leb{2}_{q,p}}^2 +s^{\beta(1-2\alpha+3\theta)} \norm{\gamma}{\Leb{\infty}_{q,p}}\norm{\nabla_q \gamma}{\Leb{\infty}_{q,p}}+ s^{\beta(2-2\alpha-\theta)}\norm{|p|^2 \gamma}{\Leb{\infty}_{q,p}}^2\\
 &\apprle \varepsilon_1^2 (s^{\beta(5-2\alpha)} + s^{\beta(1-2\alpha+3\theta)+\theta_1}+s^{\beta(2-2\alpha-\theta)+4(2-2\alpha)})\apprle \varepsilon_1^2 s^{\theta^*},
\end{aligned}
\end{equation}
where
\[ \theta^* = \min \{\beta(5-2\alpha),\beta(1-2\alpha+3\theta)+\theta_1, \beta(2-2\alpha-\theta)+8-8\alpha \}<0.\]}
Here we choose $\alpha>1$ sufficiently close to $1$ so that $\theta_1<0$ and $\theta^*<0$. 
Hence by \eqref{eq:commutator-relationship}, \eqref{eq:nabla-q-E-epsilon-2}, \eqref{eq:nabla-q-gamma-2}, we have
\begin{align*}
\norm{p^m\nabla_q \gamma(s)}{\Leb{\infty}_{q,p}}&\leq \norm{|p|\nabla_q \gamma(1)}{\Leb{\infty}_{q,p}}+ \int_s^1 \left(\norm{\nabla_q E}{\Leb{\infty}_q}\norm{|p|\nabla_p \gamma}{\Leb{\infty}_{q,p}} + \norm{E}{\Leb{\infty}_q}\norm{\nabla_q \gamma}{\Leb{\infty}_{q,p}}\right)\tau^{1-2\alpha} d\tau\\
&\apprle \varepsilon_1+\varepsilon_1^4 \int_s^1 \tau^{\theta^*+2-2\alpha+1-2\alpha} d\tau +\varepsilon_1^2 \int_s^1 \tau^{\theta_1+1-2\alpha} d\tau\\
&\apprle \varepsilon_1 (s^{\theta^* +4-4\alpha}+s^{\theta_1+2-2\alpha})\\
&\apprle \varepsilon_1 s^{\theta_2},
\end{align*}
where
\[ \theta_2=\min\{\theta^*+4-4\alpha,\theta_1+2-2\alpha\}	\]
which is still negative if we choose $\alpha>1$ sufficiently close to $1$, which implies the existence of $\delta>0$ in Proposition \ref{prop:bootstrap}. This completes the proof of Proposition \ref{prop:bootstrap}.\end{proof}

The following corollary is useful when we study the asymptotic behavior of solutions. 
\begin{corollary}\label{cor:E0-existence}
Let $\gamma$ be a solution of as in Proposition \ref{prop:bootstrap}, which is moreover defined on $(0,1]$. Then the limit
\[ E_0(q):=\lim_{s\rightarrow 0+} E(s,q) \]
exists and is bounded by
\begin{equation}\label{eq:E0-control-original}
 \norm{E_0}{\Leb{\infty}_q}\apprle \varepsilon_1^2.
\end{equation}
In addition, we have the following convergence rate: if $0\leq s_1< s_2\leq 1$, then 
\begin{equation}\label{eq:conv-rate}
 \norm{E(s_1)-E(s_2)}{\Leb{\infty}_q}\apprle \varepsilon^2 s_2^{\beta^\prime}\action{\ln s_2},
\end{equation} 
where $\beta^\prime=1{{0-9\alpha}}$ if $\alpha>1$ and $\beta^\prime=2\alpha-1$ if $\alpha<1$.  
\end{corollary}

Now we are ready to prove Theorem \ref{thm:A}.
\begin{proof}
{We only prove the case $\alpha>1$ since the case for $\alpha<1$ is simpler. We will show that there exists a sufficiently small $\varepsilon_0>0$ such that if $\mu_0$ satisfies
\[ B=\norm{\mu_0}{\Leb{2}_{x,v}}+\norm{\action{x}^4\mu_0}{\Leb{\infty}_{x,v}}+\norm{\action{v}^4\mu_0}{\Leb{\infty}_{x,v}}+\norm{\nabla_{x,v}\mu_0}{\Leb{2}_{x,v}\cap\Leb{\infty}_{x,v}}\leq \varepsilon_0,\]
then there exists a unique global solution $\mu$ to the problem. We first note that by Theorem \ref{thm:VR-LWP}, there exists $\varepsilon_1>0$ such that if $B\leq \varepsilon_1$, there exists $T>1$ such that  the problem admits a unique local solution $\mu$ satisfying 
\[ \sup_{t\in [0,T]} \norm{\mu(t)}{\Leb{2}_{x,v}}=\norm{\mu_0}{\Leb{2}_{x,v}},\quad \sup_{t\in [0,T]} (\norm{\action{x}^4\mu(t)}{\Leb{\infty}_{x,v}}+\norm{\action{v}^4\mu(t)}{\Leb{\infty}_{x,v}}+\norm{\nabla_{x,v}\mu(t)}{\Leb{2}_{x,v}\cap\Leb{\infty}_{x,v}})\leq CB.\]

If we use the pseudoconformal transform defined in \eqref{eq:pseudo-conformal-relationship}, then $\gamma_1(q,p)=\gamma(1,q,p)=\mu(1,q,q-p)$ satisfies
\begin{align*}
\norm{\gamma_1}{\Leb{2}_{q,p}}+\norm{\action{p}^4\gamma_1}{\Leb{\infty}_{q,p}}+\norm{\nabla_{p,q}\gamma_1}{\Leb{\infty}_{q,p}}\apprle B.
\end{align*}
Hence if we choose $\varepsilon_0>0$ sufficiently small, then it follows from Proposition \ref{prop:bootstrap} that the local solution $\gamma$ becomes global on $(0,1]$. Moreover, it follows from Corollary \ref{cor:E0-existence} that $E_0(q):=\lim_{s\rightarrow 0+}E(s,q)$ exists and satisfies \eqref{eq:E0-control-original} and \eqref{eq:conv-rate}. This proves that the local solution becomes global on $[0,\infty)$. }

 To study the asymptotic behavior of solutions, we may assume that {$\norm{\action{p}\nabla_{q,p}\gamma_1}{\Leb{\infty}_{q,p}}\leq \varepsilon$} since the local uniform convergence can be proved via a simple modification. Define 
\begin{align*}
 \Phi_\alpha(s,q,p)&=(\Phi^1_\alpha(s,q,p),\Phi^2_\alpha(s,q,p))=\left(q+sp+\frac{\lambda}{2-2\alpha}s^{3-2\alpha}E_0(q),p+\frac{\lambda}{2-2\alpha} s^{2-2\alpha}E_0(q)\right)
\end{align*}
and  
\[ \nu(s,q,p)=\gamma(s,\Phi_\alpha(s,q,p)).\]

A direct computation gives
\begin{align*}
\partial_s \nu(s,q,p)&=\lambda \frac{3-2\alpha}{2-2\alpha} E_0(q)\cdot(\nabla_q \gamma)(s,\Phi_\alpha(s,q,p))\\ 
&\relphantom{=}+\lambda s^{1-2\alpha}\left(E_0(q)-E\left(s,\Phi_\alpha^1(s,q,p)\right)\right)\cdot (\nabla_p \gamma)(s,\Phi_\alpha(s,q,p)).
\end{align*}

By \eqref{eq:weight-p-a}, \eqref{eq:nabla-q-E-epsilon-2}, and \eqref{eq:E0-control-original}, we have
\begin{equation}\label{eq:growth-control}
\begin{aligned}
&s^{1-2\alpha}\left|\left(E_0(q)-E\left(s,\Phi_\alpha^1(s,q,p)\right)\right)\cdot (\nabla_p \gamma)(s,\Phi_\alpha(s,q,p))\right|\\
&\apprle s^{1-2\alpha}|E_0(q)-E(s,q)|\norm{\nabla_p \gamma}{\Leb{\infty}_{q,p}}+s^{1-2\alpha}|[E(s,q)-E(s,\Phi_\alpha^1(s,q,p))]\cdot \nabla_p \gamma| \\
&\apprle \varepsilon^2 s^{1-2\alpha} {{s^{10-9\alpha} }}\action{\ln s}+\varepsilon^2 s^{1-2\alpha}\norm{\nabla_q E}{\Leb{\infty}}\left| sp\cdot \nabla_p \gamma + \frac{\lambda}{2-2\alpha} s^{3-2\alpha} E_0(q) \cdot \nabla_p \gamma\right|\\
&\apprle \varepsilon^2({{s^{11-11\alpha}}}\action{\ln s}+ s^{4-4\alpha+\theta^*})
\end{aligned}   
\end{equation}
and 
\[
\norm{E_0(q)\cdot\nabla_q \gamma}{\Leb{\infty}_q}\apprle \varepsilon^2 s^{\theta_1}.
\] 
Hence it follows that 
\begin{equation}\label{eq:convergence-pdf}
	\begin{split}
		\norm{\partial_s \nu}{\Leb{\infty}_{q,p}}\apprle \varepsilon^2 \left(s^{\theta_1}+s^{11-11\alpha}\action{\ln s}+ s^{4-4\alpha+\theta^*} \right).
	\end{split}
\end{equation}

Since $11-11\alpha>-1$, $4-4\alpha+\theta^*>-1$, and $\theta_1>-1$, the right hand side is integrable on $(0,1]$ which implies the existence of scattering.      This completes the proof of Theorem \ref{thm:A}.
\end{proof}

\begin{remark}\label{remark:threshold}\leavevmode
\begin{enumerate}
\item Write 
\[ \mathcal{A}:(s,q,p)\mapsto \left(s,q+ps+\frac{\lambda}{2-2\alpha} s^{3-2\alpha} E_0(q), p+ \frac{\lambda}{2-2\alpha} s^{2-2\alpha}E_0(q)\right) \]
and let $\mathcal{I}$ be the pseudo-conformal inversion map defined by 
\[ \mathcal{I}:(t,x,v)\mapsto \left(\frac{1}{t},\frac{x}{t},x-tv\right).\]
Since 
\[ (\gamma \circ \mathcal{A})(s,q,p)=(\mu \circ \mathcal{I}\circ \mathcal{A})(s,q,p),\]
it follows that 
\[  (\gamma\circ\mathcal{A})(s,q,p)=\mu\left(\frac{1}{s},\frac{q}{s}+p +\frac{\lambda}{2-2\alpha}s^{2-2\alpha}E_0(q),q\right).\]
Hence, by changing the coordinates, we proved that there exist $E_\infty$ and $\mu_\infty$ satisfying 
\[  \mu\left(t,x+tv+\frac{\lambda}{2-2\alpha}t^{2\alpha-2} E_\infty(v),v\right)\rightarrow \mu_\infty(x,v) \]
locally uniform in $(x,v)$ as $t\rightarrow \infty$.  {{That $\mu_{\infty}\in L^2_{x,v}$ follows from the conservation of mass and Fatou's lemma.
}}
\item {{We now collect all the conditions on the parameter $\alpha$ that have appeared throughout the paper.  
\[ \begin{cases}
\beta(1-2\alpha+3\theta)+1-2\alpha>-1, & \\
-1<\theta_1=\min\{\beta(5-2\alpha)+2-2\alpha,\beta(2-2\alpha-\theta)+10-10\alpha \}<0,  &\\
4-4\alpha+\theta^*>-1,  &\\
11-11\alpha>-1, &\\
2 - 2\alpha < \theta < \tfrac{2\alpha - 1}{3} .
\end{cases}
\]
where
\[ \theta^* = \min \{\beta(5-2\alpha),\beta(1-2\alpha+3\theta)+\theta_1, \beta(2-2\alpha-\theta)+8-8\alpha \}.\]
}

Our goal is to choose $\beta < 0$ and $2 - 2\alpha < \theta < \tfrac{2\alpha - 1}{3}$ so that the above system admits a solution for all $1<\alpha<1+\delta$. One can easily see that for $\delta\ll 1$, this linear system must have a solution since $\alpha=1, \theta=1/10$ is a solution in extreme case. Indeed, with the help of Mathematica, if we choose $\beta = -0.135$, then we can choose $\delta = 0.0801$. }

\item From \eqref{eq:growth-control}, even we assume the best scenario, say $\norm{\nabla_q E}{\Leb{\infty}}$, $\norm{\nabla_p \gamma}{\Leb{\infty}}$ are bounded in $s$, we have a growth $s^{4-4\alpha}$ which is not integrable over $(0,1]$ if  $\alpha \geq 5/4$. It would be interesting whether we can prove the existence of modified scattering when $1<\alpha<5/4$ in an appropriate topology. 
\end{enumerate}
\end{remark}

\section{Existence of wave operators}\label{sec:4}
In this section, we construct a wave operator for the problem \eqref{eq:VR}.  Observe that the system \eqref{eq:VR-pct} can be written as 
\begin{equation}\label{eq:VR-pct-PB}
 \gamma_s + \{\gamma,\mathcal{H}\}=0,\quad \{f,g\}:=\nabla_q f \cdot \nabla_p g - \nabla_p f \cdot \nabla_q g.
\end{equation}
with the Hamiltonian
\[ \mathcal{H}(s,q,p):=\frac{|p|^2}{2}-\lambda s^{1-2\alpha} \phi(s,q),\quad \phi(s,q)=c_\alpha \iint_{\mathbb{R}^3_y\times\mathbb{R}^3_u} \frac{\gamma^2(s,y,u)}{|q-y|^{5-2\alpha}}dydu.\]
Note that the Hamiltonian $\mathcal{H}$ contains a strong singularity at $s=0$. 

To solve \eqref{eq:VR-pct-PB} at ``$s=0$'', we first find a new coordinate system that preserves the symplectic structure and mitigates the singularity at $s=0$. It is obtained by the type-3 generating function (see e.g. \cite[Sections 3 and 4]{C01}):
\[ S(s,w,p):=w\cdot p+\frac{|p|^2}{2}s-\lambda\frac{s^{2-2\alpha}}{2-2\alpha}\phi_0(w), \]
where $\phi_0(w)=\phi(0,w)$. It gives us the following change of variables:
\begin{align}\label{COV}
	&q=w+zs+\lambda \frac{s^{3-2\alpha}}{2-2\alpha}\nabla\phi_0(w),&&w=q-ps \\
	&p=z+\lambda\frac{s^{2-2\alpha}}{2-2\alpha}\nabla\phi_0(w),&&z=p-\lambda\frac{s^{2-2\alpha}}{2-2\alpha}\nabla\phi_0(q-ps)\label{COV-2}
\end{align}
{{Note also that in the original coordinates $(t,x,v)$, we have
\begin{equation}
	\begin{split}
z=x-tv-\lambda\frac{t^{2-2\alpha}}{2-2\alpha}E_0(v),\quad w=v.		
	\end{split}
\end{equation} which is consistent with the asymptotic behavior. } }

We can compute 
\begin{equation}\label{JacobChangeCoord}
\frac{\partial (w,z)}{\partial (q,p)}=\begin{pmatrix} \mathrm{Id}&-s\mathrm{Id}\\-\lambda\frac{s^{2-2\alpha}}{2-2\alpha}\nabla E_0&\mathrm{Id}+\lambda \frac{s^{3-2\alpha}}{2-2\alpha}\nabla E_0\end{pmatrix},
\end{equation}
where $E=\nabla \phi$ and $E_0=\nabla \phi_0$. 
This corresponds to the new Hamiltonian
\begin{align*}
\mathcal K(s,w,z) &:= \mathcal H(s,q,p) - \partial_sS(s,w,p) = -\lambda s^{1-2\alpha} \big[ \phi(s,q) -  \phi_0(w)\big].
\end{align*}
which has milder singularity at ``$s=0$'' than $\mathcal{H}$. Also, the vector fields are 
\begin{equation}\label{eq:dwK}
\begin{split}
\nabla_w\mathcal{K}&=-\lambda s^{1-2\alpha}\left\{E(s,q)-E_0(w)\right\}-\lambda^2\frac{s^{4-4\alpha}}{2-2\alpha} E(s,q)\cdot\nabla E_0(w),\\\qquad\nabla_z\mathcal{K}&=-\lambda s^{2-2\alpha}E(s,q),
\end{split}
\end{equation}
Hence if we define
\begin{equation}\label{eq:sigma-gamma}
\sigma(s,w,z):=\gamma(s,q,p), 
\end{equation}
then $\sigma$ satisfies
\begin{equation}\label{EqSigma}
	0=\partial_s\sigma+\{\sigma,\mathcal{K} \}=\partial_s\sigma+\nabla_w\sigma\cdot\nabla_z\mathcal{K}-\nabla_z\sigma\cdot\nabla_w\mathcal{K}.
\end{equation}

To state our theorem, we define an electric field 
\begin{align*}
E[\sigma](s,Q)&=E(s,Q)=c_\alpha'\iint_{\mathbb{R}^3_w\times\mathbb{R}^3_z} \frac{Q-q(s,w,z)}{|Q-q(s,w,z)|^{5-2\alpha}} \sigma^2(s,w,z)dwdz\\
&=c_\alpha'\iint_{\mathbb{R}^3_q\times\mathbb{R}^3_p} \frac{Q-q}{|Q-q|^{5-2\alpha}} \gamma^2(s,q,p)dqdp.
\end{align*}
where we used a change of variable \eqref{eq:sigma-gamma}. One can easily see that for $0<s_1<s_2$ and $Q\in\mathbb{R}^3$, we have 
\begin{align*}
|E(s_2,Q)-E(s_1,Q)|&\apprle R^{2\alpha-1} \norm{\action{z}^4\sigma}{\Leb{\infty}_{s,w,z}}^2+R^{2\alpha-4} \norm{\sigma}{\Leb{\infty}_s \Leb{2}_{w,z}}\norm{\sigma(s_2)-\sigma(s_1)}{\Leb{2}_{w,z}}.
\end{align*}
uniformly in $R>0$. Since $\alpha>1/2$, we have
\begin{align*}
|E(s_2,Q)-E(s_1,Q)|&\apprle \norm{\action{z}^4\sigma}{\Leb{\infty}_{s,w,z}}^{(8-4\alpha)/3} (\norm{\sigma}{\Leb{\infty}_s \Leb{2}_{w,z}}\norm{\sigma(s_2)-\sigma(s_1)}{\Leb{2}_{w,z}})^{(2\alpha-1)/3}.
\end{align*}
which implies that $E$ is continuous in $s$ if $\sigma$ is continuous in $s$.

Now we are ready to present the main theorem of this section.
\begin{theorem}\label{thm:Wave}
{{Let {$\alpha \in (1/2,35/33)\setminus\{1\}$}.}} Assume that initial data $\sigma_0$ and {{$E_0(q):=E[\sigma](0,q)$}} satisfy 
\begin{equation}\label{eq:E0-control}
\norm{E_0}{\Sob{3}{\infty}}\leq c_0^2,
\end{equation}
and 
\begin{equation}\label{eq:sigma0-control}
\norm{\sigma_0}{\Leb{2}_{w,z}}+{\norm{\langle z\rangle\nabla_{z,w}\sigma_0}{\Leb{2}_{w,z}}}+\norm{\action{z}^5\sigma_0}{\Leb{\infty}_{w,z}}+\sum_{0\leq m+n\leq 2} \norm{\action{z}^m \nabla_z^m \nabla_w^n \sigma_0}{\Leb{\infty}_{w,z}}\leq c_0.
\end{equation}
Then there exist $T^*=T^*(c_0)>0$ and a unique strong solution $\sigma \in C([0,T^*];\Leb{2}_{w,z})$ of \eqref{EqSigma} with initial data $\sigma(s=0)=\sigma_0$ such that $s\partial_s \sigma$, $\nabla_{w,z}\sigma \in C_{s,w,z}$. Moreover, for $0\leq s<T^*$, we have that for any $\ell \in \mathbb{N}$,
\begin{align}
&\norm{\sigma(s)}{\Leb{2}_{w,z}}+\norm{\action{z}^5 \sigma(s)}{\Leb{\infty}_{w,z}}+\norm{\nabla_{w,z}\sigma(s)}{\Leb{\infty}_{w,z}}+{ \norm{\nabla_{w,z}\sigma(s)}{\Leb{2}_{w,z}}}\apprle c_0,\label{eq:sigma-c0} \\
&{\norm{\action{w,z}^l \sigma(s)}{\Leb{r}_{w,z}}\apprle \norm{\action{w,z}^l \sigma_0}{\Leb{r}_{w,z}},\quad r\in \{2,\infty\}.}\label{eq:sigma-sigma0}
\end{align}
{ Moreover, if $c_0$ is sufficiently small, we may take $T^*=1$.}
\end{theorem}
The solution will be constructed by a Picard iteration. To do this, we need to get a priori estimates of particle density functions, which leads us to get a uniform estimate for sequences of functions. We assume  the following bootstrap assumptions on $\sigma$ to get uniform estimates of $\sigma$ and the derivatives of $\mathcal{K}$:
\begin{subequations}\label{eq:wave-bootstrap}
\begin{align}
\norm{\sigma(s)}{\Leb{2}_{w,z}}+\norm{\action{z}^5 \sigma(s)}{\Leb{\infty}_{w,z}}\leq {A_1}\leq 4c_0,\label{eq:wave-bootstrap-1}\\
{\norm{\nabla_w\sigma(s)}{\Leb{2}_{w,z}\cap\Leb{\infty}_{w,z}}+\norm{\theta\nabla_z \sigma(s)}{\Leb{2}_{w,z}\cap\Leb{\infty}_{w,z}}}\leq {A_2}\leq 4c_0,\label{eq:wave-bootstrap-2} \\
\norm{\nabla_{w,w}^2\sigma(s)}{\Leb{\infty}_{w,z}}+\norm{\theta\nabla_{w,z}^2\sigma(s)}{\Leb{\infty}_{w,z}}+\norm{\theta^2\nabla_{z,z}^2\sigma(s)}{\Leb{\infty}_{w,z}}\leq {A_3}\leq 4c_0,\label{eq:wave-bootstrap-3}
\end{align}
\end{subequations}
where $\theta$ is defined by
\begin{align*}
	\theta(s,z):=\frac{\langle z\rangle}{1+s\langle z\rangle}, \quad \frac{1}{2}\min\{\langle z\rangle,s^{-1} \}\leq \theta(s,z)\leq \min\{\langle z\rangle,s^{-1} \}
\end{align*}
which satisfies the following differential equations:
\begin{align*}
	\partial_s\theta=-\theta^2\quad\text{and}\quad \nabla_z\theta=\left(\frac{z}{\langle z\rangle^3} \right)\cdot \theta^2.
\end{align*}

In Subsection \ref{subs:comm-relationship}, we summarize the commutation relations with the operator $\mathfrak{L}=\partial_s + \{\cdot,\mathcal{K}\}$. Using this relation, we prove a priori estimates on the electric field and particle density function $\sigma$ in Subsections \ref{subs:electric-field} and \ref{subs:particle-density}. The proof of Theorem \ref{thm:Wave} will be given in Subsection \ref{subs:wave-theorem}. Finally, we prove Theorem \ref{thm:B} in Subsection \ref{subsec:propagation}.

\subsection{Commutation relations}\label{subs:comm-relationship}
Writing $\mathfrak{L}=\partial_s+\{\cdot,\mathcal{K}\}$, for the moments in $w,z$, we have the commutation relations
\begin{equation}\label{eq:commutation-relations-VR-wave}
\begin{aligned}
	&\mathfrak{L}[w_j\sigma]=-\lambda s^{2-2\alpha}E_{j}(q){{\sigma}} ,\\
	&\mathfrak{L}[{{z_j}}\sigma]={{\lambda s^{1-2\alpha}\left\{E_{j}(s,q)-E_{0,j}(w)\right\}+\lambda^2\frac{s^{4-4\alpha}}{2-2\alpha} E(s,q)\cdot\nabla E_{0,j}(w)}}.
\end{aligned}
\end{equation}
For the first-order differential operator $\partial_1,\partial_2$ in $w$ or $z$, we have 
\begin{equation}\label{eq:GgradComm}
\begin{split}
\mathfrak{L}(\partial_1\sigma)&=\{\partial_1\mathcal{K},\sigma\},\\
 \mathfrak{L}(\partial_2\partial_1\sigma)&=\{\partial_1\mathcal{K},\partial_2\sigma\}+\{\partial_2\mathcal{K},\partial_1\sigma\}+\{\partial_2\partial_1\mathcal{K},\sigma\}.
\end{split}
\end{equation}

In other words, we can write 
\begin{equation}\label{eq:gradComm}
\mathfrak L \begin{pmatrix} \nabla_{w} \sigma \\ \nabla_{z} \sigma  \end{pmatrix} = \begin{pmatrix} -\nabla_{w}\nabla_{z} \mathcal K & \nabla_{w}^2 \mathcal K \\
-\nabla_{z}^2 \mathcal K&  \nabla_{w}\nabla_{z} \mathcal K  \end{pmatrix}\begin{pmatrix} \nabla_{w} \sigma \\ \nabla_{z} \sigma   \end{pmatrix}, 
\end{equation}
with
\begin{equation}\label{FormulasHessianK}
\begin{aligned}
\nabla_{w^jw^k}^2 \mathcal K &= -\lambda s^{1-2\alpha}[\partial_j E_k(s,q) -\partial_j E_{0,k}(w)]\\
&\relphantom{=} -\lambda^2\frac{s^{4-4\alpha}}{2-2\alpha}\left\{\partial_j E(s,q)\cdot\partial_k E_{0}(w) \right.\\
&\relphantom{=} \left.+\partial_kE(s,q)\cdot\partial_j E_{0}(w) +\partial_j\partial_k E_0(w)\cdot E(s,q)\right\}\\
&\relphantom{=} -\lambda^3 \frac{s^{7-6\alpha}}{(2-2\alpha)^2} \partial_k E_{0,a}(w)\partial_jE_{0,b}(w) \cdot \partial_a E_b(s,q),\\
\nabla_{w}\nabla_{z} \mathcal K &=-\lambda s^{2-2\alpha} \nabla E(s,q)  - \lambda^2 \frac{s^{5-4\alpha}}{2-2\alpha} (\nabla E(s,q)\cdot \nabla) E_0(w), \\
\nabla_{z,z}^2 \mathcal K &= - \lambda s^{3-2\alpha} \nabla E(s,q),
\end{aligned}
\end{equation}
where we used Einstein summation convention for repeated indices.

\subsection{A priori estimates on electric field}\label{subs:electric-field}
In this subsection, we obtain a priori estimates on the electric field under the bootstrap assumption \eqref{eq:wave-bootstrap}.

We first consider the case $\alpha>1$ since the case for $\alpha<1$ is similar, and in fact, it is simpler.
\begin{lemma}\label{lem:modified-a-priori-electric-field}
{{Let $1<\alpha<23/20$}} and $\sigma \in C([0,T];\Leb{2}_{w,z})$ with $\sigma(0)=\sigma_0$ satisfy {
\begin{equation}\label{eq:continuity-eqn-sigma}
   \partial_s \rho + \Div_q \mathbf{j}=0\quad \text{in } [0,T]\times \mathbb{R}^3,
\end{equation}
where 
\[ \rho(s,q)=\int_{\mathbb{R}^3} \gamma^2(s,q,p)dp\quad \text{and}\quad \boldj(s,q)=\int_{\mathbb{R}^3} p \gamma^2 (s,q,p)dp. \]
}
Suppose that $E_0$ satisfies \eqref{eq:E0-control}. Then there exists $T^*(c_0)\in (0,T]$ such that 
\begin{enumerate}[label=\textnormal{(\roman*)}]
\item Assume that \eqref{eq:wave-bootstrap-1} holds. Then 
\begin{equation}\label{eq:modified-K-control-1}
\begin{aligned}
|\nabla_z \mathcal{K}(w,z)|&\apprle  c_0^2s^{2-2\alpha},\\
|\nabla_w \mathcal{K}(w,z)|&\apprle  c_0^2 \{s^{2-2\alpha} \min\{s^{-1},|z|\}+s^{5-5\alpha} \} 
\end{aligned}
\end{equation}
for $0< s\leq T^*$.
\item Assume that \eqref{eq:wave-bootstrap-1} and \eqref{eq:wave-bootstrap-2} hold. Then 
\begin{equation}\label{eq:modified-K-control-2}
\norm{\theta\nabla^2_{z,z} \mathcal{K}}{\Leb{\infty}_{w,z}}+\norm{\nabla_{z,w}^2 \mathcal{K}}{\Leb{\infty}_{w,z}}+\norm{\theta^{-1} \nabla^2_{w,w} \mathcal{K}}{\Leb{\infty}_{z,w}}\apprle c_0^2 s^{\frac{20}{3}(1-\alpha)}
\end{equation}
for $0< s\leq T^*$.
\end{enumerate} 
\end{lemma}
\begin{proof}
(i) Recall from \eqref{COV-2} and  \eqref{JacobChangeCoord} that 
\begin{equation}\label{eq:modified-p-z-q-w}
	{{|p-z|\lesssim c_0^2}} s^{2-2\alpha}, \quad |q-w|\leq s|z|+\frac{1}{2-2\alpha}s^{3-2\alpha}c_0^2.
\end{equation}
hold and the change of variable preserves volume. These imply that 
\begin{equation}\label{eq:change-of-variable-sigma-gamma}
	\begin{split}
		&\norm{\gamma(s)}{L^r_{q,p}}=\norm{\sigma(s)}{L^r_{q,p}}, \\
		&\norm{|p|^{\beta}\gamma(s)}{L^r_{q,p}}\apprle {{ c_0^{2\beta}}} s^{(2-2\alpha)\beta}\norm{\sigma(s)}{L^r_{w,z}}+\norm{|z|^{\beta}\sigma(s)}{L^r_{w,z}}.
	\end{split}	
\end{equation}

In addition, since $\frac{\partial z}{\partial p}=\mathrm{Id}+O(c_0 s^{3-2\alpha})$ has bounded Jacobian, we see that 
\begin{equation}\label{eq:modified-j-estimate}
	\norm{\bold{j}(s)}{L^{\infty}_q}\leq \norm*{\int [|z|+c_0^2 s^{2-2\alpha}]\sigma^2 dz}{L^{\infty}_w}\apprle {{\action{c_0}^2s^{2-2\alpha}}}\norm{\langle z\rangle^3\sigma}{L^{\infty}_{w,z}}^2
\end{equation}
and {{using \eqref{continuityE2}}}, we obtain that for $2^{-k-1}\leq s_2\leq s_1\leq 2^{-k}$, 
\begin{equation*}
	\norm{E(s_2,q)-E(s_1,q)}{L^{\infty}_q}\apprle  \langle c_0\rangle^2 2^{-k(4-3\alpha)}(\norm{\langle z\rangle^3\sigma}{L^{\infty}_{s,w,z}}^2+\norm{\sigma}{L^{\infty}_sL^2_{w,z}}^2)
\end{equation*}
Since $\alpha<4/3$, we see that $E(2^{-k})$ is Cauchy in $L^{\infty}_q$ and  hence it follows that
\begin{equation}\label{modified-convergeE}
	\norm{E(s,q)-E_0(q)}{L^{\infty}_q}\apprle  \langle c_0\rangle^2 {{{s}^{2-2\alpha} }}s^{2-\alpha}(\norm{\langle z\rangle^3\sigma}{L^{\infty}_{s,w,z}}^2+\norm{\sigma}{L^{\infty}_sL^2_{w,z}}^2)
\end{equation}
so that we have 
\begin{equation}\label{eq:modified-energy-control}
	|E(s,q)|\apprle  |E(s,q)-E_0(q)|+|E_0(q)|\apprle  c_0^2
\end{equation}
Then by \eqref{eq:dwK}, \eqref{eq:modified-p-z-q-w}, \eqref{modified-convergeE}, and \eqref{eq:modified-energy-control}, we have
\begin{equation*}
\begin{split}
&|\nabla_z \mathcal{K}(w,z)|\apprle  c_0^2 s^{2-2\alpha}, \\
&|\nabla_w \mathcal{K}(w,z)|=\left|-\lambda s^{1-2\alpha}\{E(s,q)-E_0(q) \}-\lambda s^{1-2\alpha}\{E_0(q)-E_0(w) \}+O(c_0^4 s^{4-4\alpha})\right| \\
&\apprle c_0^2  s^{2-2\alpha}\min\{|z|,s^{-1}\} +O(c_0^2s^{5-5\alpha}).
\end{split}
\end{equation*}

(ii)  By \eqref{COV}, we see that
\begin{equation*}
	\begin{split}
		\nabla_q\gamma=\nabla_w\sigma-\lambda \frac{s^{2-2\alpha}}{2-2\alpha}\nabla E_0\cdot \nabla_z\sigma.
	\end{split}
\end{equation*}
so that for $r \in \{2,\infty\}$, it follows from \eqref{eq:E0-control} that 
\begin{equation*}
	\begin{split}
		\norm{\nabla_q\gamma}{L^r_{q,p}}\apprle  c_0^2 {{s^{2-2\alpha}}}\norm{\nabla_{w,z}\sigma}{L^r_{w,z}}
	\end{split}
\end{equation*}
and 
\begin{equation*}
	\begin{split}
		\norm{\nabla_q\bold{j}(s)}{L^{\infty}_q}&\apprle  c_0^2 {{s^{2-2\alpha}}} \left\Vert \int [|z|+c_0^2s^{2-2\alpha}]|\sigma|\cdot|\nabla_{w,z}\sigma | dz  \right\Vert_{L^{\infty}_q} \\
		&\apprle {{ c_0^2}} {{s^{4-4\alpha}}}(\norm{\langle z\rangle^5\sigma}{L^{\infty}_{z,w}}^2+\norm{\nabla_{w,z}\sigma}{L^{\infty}_{w,z}}^2 )
	\end{split}
\end{equation*}
Applying the same argument as before by \eqref{continuityDE2}, we can get
		{{\begin{align*}
			\norm{\nabla_qE(s,q)-\nabla_q E_0(q)}{L^{\infty}_q}&\apprle  c_0^2 s^{\frac{17-14\alpha}{3}}\left(\norm{\nabla_{w,z}\sigma}{\Leb{\infty}_{s,w,z}}^2+\norm{\action{z}^5\sigma}{\Leb{\infty}_{s,w,z}}^2+\norm{\sigma}{\Leb{\infty}_{s}\Leb{2}_{w,z}}^2 \right) \\
			\norm{\nabla E(s,\cdot)}{L^{\infty}_q}&\lesssim c_0^2
		\end{align*}}}
by choosing $T(c_0)>0$ sufficiently small.  Using the formulas \eqref{FormulasHessianK}, we directly see that 
	\begin{align*}
	\theta|\nabla^2_{z,z}\mathcal{K} | &\leq |\nabla E |\cdot 	s^{3-2\alpha} \min \{s^{-1},{{\action{z}}} \}\apprle c_0^2 s^{2-2\alpha} \\
	|\nabla^2_{w,z}\mathcal{K}|&{{\lesssim }}s^{2-2\alpha} |\nabla E|(1+s^{3-2\alpha}|\nabla E_0 | )\apprle 2c_0^2 s^{2-2\alpha}.
		\end{align*}
Moreover, it follows from  \eqref{FormulasHessianK} that by choosing $T(c_0)>0$ sufficiently small, we have
	{{\begin{align*}
	|\nabla_{w,w}^2\mathcal{K} |&\apprle s^{1-2\alpha}|\nabla E_0(q)-\nabla E_0(w) |+ s^{1-2\alpha}|\nabla E(s,q)- \nabla E_0(w) | \\
	&+s^{4-4\alpha}[{2|\nabla E_0 |}|\nabla E|+|\nabla^2 E_0 | |E|   ]+s^{7-6\alpha}|\nabla E_0 |^2|\nabla E | \\
	&\lesssim c_0^2  s^{2-2\alpha}\min\{s^{-1},|z|\}+c_0^2 s^{4-4\alpha}+c_0^2 s^{\frac{20}{3}(1-\alpha)}+s^{7-6\alpha}c_0^2 .
	\end{align*}}}
which implies \eqref{eq:modified-K-control-2}. This completes the proof of Lemma \ref{lem:modified-a-priori-electric-field}.
\end{proof}

Following a similar argument, we can also prove the following lemma for the case $\alpha<1$.

\begin{lemma}\label{lem:a-priori-electric-field}
Let $1/2<\alpha<1$ and let $\sigma \in C([0,T];\Leb{2}_{w,z})$ with $\sigma(0)=\sigma_0$ satisfy \eqref{eq:continuity-eqn-sigma}. Suppose that $E_0$ satisfies \eqref{eq:E0-control}. Then there exists $T^*(c_0)\in (0,T]$ such that 
\begin{enumerate}[label=\textnormal{(\roman*)}]
\item Assume that \eqref{eq:wave-bootstrap-1} holds. Then 
\begin{equation}\label{eq:K-control-1}
\begin{aligned}
&|\nabla_z \mathcal{K}(w,z)|\apprle c_0^2,\\
&|\nabla_w \mathcal{K}(w,z)|\apprle c_0^2 \{s^{2-2\alpha} {{\min\{s^{-1},|z|\}+1}}\} \leq c_0^2 \min\{s^{-1},\action{z}\}.
\end{aligned}
\end{equation}
for $0< s\leq T^*$.
\item Assume that \eqref{eq:wave-bootstrap-1} and \eqref{eq:wave-bootstrap-2} hold. Then 
\begin{equation}\label{eq:K-control-2}
\norm{\theta\nabla^2_{z,z} \mathcal{K}}{\Leb{\infty}_{w,z}}+\norm{\nabla_{z,w}^2 \mathcal{K}}{\Leb{\infty}_{w,z}}+\norm{\theta^{-1} \nabla^2_{w,w} \mathcal{K}}{\Leb{\infty}_{z,w}}\apprle c_0^2.
\end{equation}
for $0< s\leq T^*$.
\end{enumerate} 
\end{lemma}

\subsection{A priori estimates on particle density}\label{subs:particle-density}
Next, we obtain a priori estimate on particle densities which is a solution of some Hamiltonian PDE that is relevant to our system.

\begin{lemma}\label{lem:particle-density-bootstrap}
{{Let $\alpha \in (1/2,35/33)\setminus\{1\}$.}} Suppose that $\sigma \in C([0,T];\Leb{2}_{w,z})$ satisfies \eqref{EqSigma} for some Hamiltonian $\mathcal{K}$ satisfying 
\begin{enumerate}[label=\textnormal{(\roman*)}]
\item  \eqref{eq:modified-K-control-1} and \eqref{eq:modified-K-control-2} if $\alpha>1$;
\item \eqref{eq:K-control-1} and \eqref{eq:K-control-2} if $\alpha<1$.
\end{enumerate}
If $\sigma_0$ satisfies \eqref{eq:E0-control} and \eqref{eq:sigma0-control}, then there exists $T(c_0)>0$ such that \eqref{eq:wave-bootstrap} hold for $A_1=A_2=A_3=2c_0$.
\end{lemma}
\begin{proof}
	We first close the bootstrap for {$A_1$}, then for {$A_1, A_2$}, and finally for {$A_1,A_2,A_3$.} We only consider the case $\alpha>1$ since the case $\alpha<1$ is simpler. The control follows from the commutation relations:
\begin{equation}\label{CommutationRelationsSigma}
\begin{split}
\mathfrak{L}(\langle z\rangle^m\sigma)&=\sigma\{\langle z\rangle^m,\mathcal{K}\}=-m \sigma\langle z\rangle^{m-2} z\cdot \nabla_w\mathcal{K} ,\\
\mathfrak{L}(\theta\nabla_z\sigma)&=(\mathfrak{L}\ln\theta)\cdot\theta \nabla_z\sigma+\theta\{\nabla_z\mathcal{K},\sigma\}\\
&=(\mathfrak{L}\ln\theta)\cdot\theta \nabla_z\sigma+\theta\nabla_z\sigma\cdot \nabla_w\nabla_z\mathcal{K}-\nabla_w\sigma\cdot \theta\nabla_z\nabla_z\mathcal{K},\\
\mathfrak{L}(\nabla_w\sigma)&=\{\nabla_w\mathcal{K},\sigma\}=\theta\nabla_z\sigma\cdot \theta^{-1}\nabla_w\nabla_w\mathcal{K}-\nabla_w\sigma\cdot\nabla_z\nabla_w\mathcal{K}.
\end{split}
\end{equation}

Then for $r\in \{2,\infty\}$, we have
\begin{equation*}
\Vert \langle z\rangle^m\sigma(s)\Vert_{L^r_{w,z}}\le \Vert \langle z\rangle^m\sigma_0\Vert_{L^r_{w,z}}+m\int_0^s\Vert \langle z\rangle^{-1}\nabla_{w}\mathcal{K}(\tau)\Vert_{L^\infty_{w,z}}\Vert \langle z\rangle^m\sigma(\tau)\Vert_{L^r_{w,z}}d\tau.
\end{equation*}
Since 
\[ \norm{\action{z}^{-1} \nabla_w \mathcal{K}(s)}{\Leb{\infty}_{w,z}}\apprle {{c_0^2s^{5-5\alpha}}} \]
and $\alpha<6/5$, we can easily propagate \eqref{eq:wave-bootstrap-1} {{for short time.
}}
\medskip

For the derivatives, we also need to control $\theta$. Note that 
\begin{equation*}
\mathfrak{L}(\ln\theta)=-\left(1+\frac{z}{\langle z\rangle^3}\nabla_w\mathcal{K}\right)\theta\apprle  {{c_0^2 s^{5-5\alpha}.}}
\end{equation*}

Then it follows from \eqref{eq:modified-K-control-2} and  \eqref{CommutationRelationsSigma} that
\begin{equation*}
\begin{split}
\Vert \theta\nabla_z\sigma(s)\Vert_{L^r_{w,z}}&\le\Vert \theta\nabla_z\sigma_0\Vert_{L^r_{w,z}}+ c_0^2\int_0^s \tau^{\frac{20}{3}(1-\alpha)}\left\{\Vert \theta\nabla_z\sigma(\tau)\Vert_{L^r_{w,z}}+\Vert\nabla_w\sigma(\tau)\Vert_{L^r_{w,z}}\right\}d\tau,\\
\Vert \nabla_w\sigma(s)\Vert_{L^r_{w,z}}&\le\Vert \nabla_w\sigma_0\Vert_{L^r_{w,z}}+ c_0^2\int_0^s\tau^{\frac{20}{3}(1-\alpha)}\left\{\Vert \theta\nabla_z\sigma(\tau)\Vert_{L^r_{w,z}}+\Vert\nabla_w\sigma(\tau)\Vert_{L^r_{w,z}}\right\}d\tau.\\
\end{split}
\end{equation*}
Since $\alpha<23/20$, these allow us to propagate \eqref{eq:wave-bootstrap-2} for short time.

\medskip

We now propagate higher order derivatives to bound the bootstrap for $A_3$. First by interpolation in \eqref{eq:wave-bootstrap}, we observe that
\begin{equation*}
\Vert \langle z\rangle^{2.1}\nabla_{w,z}\sigma\Vert_{L^\infty_{w,z}}\le A_1+A_3,
\end{equation*}
see Appendix \ref{app:interpolation} for the proof. We will use the weight $\theta$ to control the $\partial_z$ derivatives of $\sigma$. Using \eqref{eq:GgradComm}, we find that
\begin{align*}
\mathfrak{L}(\partial_{w^j}\partial_{w^k}\sigma)&=\theta^{-1}\nabla_w\partial_{w^j}\mathcal{K}\cdot (\theta\nabla_z\partial_{w^k}\sigma)+\theta^{-1}\nabla_w\partial_{w^k}\mathcal{K}\cdot(\theta\nabla_z\partial_{w^j}\sigma)\\
&\relphantom{=}-\nabla_z\partial_{w^j}\mathcal{K}\cdot\nabla_w\partial_{w^k}\sigma-\nabla_z\partial_{w^k}\mathcal{K}\cdot\nabla_w\partial_{w^j}\sigma\\
&\relphantom{=}+\theta^{-1}\nabla_w\partial_{w^j}\partial_{w^k}\mathcal{K}\cdot(\theta\nabla_z\sigma)-\nabla_z\partial_{w^j}\partial_{w^k}\mathcal{K}\cdot\nabla_w\sigma,\\
\mathfrak{L}(\theta\partial_{z^j}\partial_{w^k}\sigma)&=\mathfrak{L}(\ln\theta)\cdot\theta\partial_{z^j}\partial_{w^k}\sigma+\theta^{-1}\nabla_w\partial_{w^k}\mathcal{K}\cdot(\theta^2\nabla_z\partial_{z^j}\sigma)\\
&\relphantom{=}-\nabla_z\partial_{w^k}\mathcal{K}\cdot(\theta\nabla_w\partial_{z^j}\sigma)+\nabla_w\partial_{z^j}\mathcal{K}\cdot (\theta\nabla_z\partial_{w^k}\sigma)\\
&\relphantom{=}-(\theta\nabla_z\partial_{z^j}\mathcal{K})\cdot\nabla_w\partial_{w^k}\sigma+\nabla_w\partial_{z^j}\partial_{w^k}\mathcal{K}\cdot(\theta\nabla_z\sigma)\\
&\relphantom{=}-\theta\nabla_z\partial_{z^j}\partial_{w^k}\mathcal{K}\cdot\nabla_w\sigma,\\
\mathfrak{L}(\theta^2\partial_{z^j}\partial_{z^k}\sigma)&=2\mathfrak{L}(\ln\theta)\cdot\theta^2\partial_{z^j}\partial_{z^k}\sigma+\nabla_w\partial_{z^k}\mathcal{K}\cdot(\theta^2\nabla_z\partial_{z^j}\sigma)\\
&\relphantom{=}-\theta\nabla_z\partial_{z^k}\mathcal{K}\cdot(\theta\nabla_w\partial_{z^j}\sigma)+\nabla_w\partial_{z^j}\mathcal{K}\cdot (\theta^2\nabla_z\partial_{z^k}\sigma)\\
&\relphantom{=}-(\theta\nabla_z\partial_{z^j}\mathcal{K})\cdot (\theta\nabla_w\partial_{z^k}\sigma)+\theta\nabla_w\partial_{z^j}\partial_{z^k}\mathcal{K}\cdot(\theta\nabla_z\sigma)\\
&\relphantom{=}-\theta^2\nabla_z\partial_{z^j}\partial_{z^k}\mathcal{K}\cdot\nabla_w\sigma,
\end{align*}
and we can proceed as for the case of one derivative once we control the new terms
\begin{equation}\label{3DerK}
\begin{split}
\Vert \theta^{-1}\nabla^3_{w,w,w}\mathcal{K}\Vert_{L^\infty_{w,z}}+\Vert \nabla^3_{w,w,z}\mathcal{K}\Vert_{L^\infty_{w,z}}+\Vert \theta\nabla^3_{w,z,z}\mathcal{K}\Vert_{L^\infty_{w,z}}+\Vert \theta^2\nabla^3_{z,z,z}\mathcal{K}\Vert_{L^\infty_{w,z}}&\apprle {{c_0^2 s^{\frac{33}{2}(1-\alpha)}}}.
\end{split}
\end{equation}

Starting from
\begin{equation*}
\begin{split}
\nabla_z\mathcal{K}=-{{\lambda s^{2-2\alpha}}} E(q),\qquad \frac{\partial q^k}{\partial z^j}=s\delta_j^k,\qquad\frac{\partial q^k}{\partial w^j}=\delta_j^k+\lambda \frac{s^{3-2\alpha}}{2-2\alpha}\partial_j\partial_k\phi_0(w),
\end{split}
\end{equation*}
where $\delta_j^k=1$ if $j=k$ and $0$ otherwise, chain rule gives
\begin{equation}\label{eq:K-third-1}
\begin{split}
\theta^2\vert \nabla^3_{z,z,z}\mathcal{K}\vert&={{s^{2-2\alpha}}} (s\theta)^2\vert\nabla^2E(q)\vert,\\
\theta\vert\nabla^3_{w,z,z}\mathcal{K}\vert&{{\lesssim s^{2-2\alpha}}} (s\theta)\cdot\left[1+s^{3-2\alpha}\vert\nabla E_0\vert\right]\cdot \vert\nabla^2E(q)\vert,\\
\vert\nabla^3_{w,w,z}\mathcal{K}\vert&\le {{s^{2-2\alpha}}}  \left[1+s^{3-2\alpha}\vert\nabla E_0\vert\right]^2\cdot \vert\nabla^2E(q)\vert\\
&\relphantom{=}+[s^{5-4\alpha}\vert\nabla E_0\vert\cdot\vert\nabla^2 E_0\vert\cdot \vert\nabla E(q)\vert,\\
\end{split}
\end{equation}
and finally, from \eqref{FormulasHessianK}, we obtain that
\begin{equation}\label{eq:K-third-2}
\begin{split}
\vert\nabla^3_{w,w,w}\mathcal{K}\vert&\le \left[s^{1-2\alpha}+ s^{4-4\alpha}\cdot \vert\nabla E_0\vert\right]\cdot \vert\nabla^2E(s,q)-\nabla^2E_0(w)\vert\\
&\relphantom{=}+s^{4-4\alpha}\cdot \left[\vert\nabla^2E\vert\cdot\vert\nabla E_0\vert+\vert\nabla E\vert\cdot \vert\nabla^2E_0\vert+\vert\nabla^3E_0\vert\cdot\vert E\vert\right]\\
&\relphantom{=}+s^{7-6\alpha}\cdot\left[\vert\nabla^2 E\vert\cdot\vert\nabla E_0\vert^2+\vert\nabla E\vert \cdot\vert\nabla E_0\vert\cdot\vert\nabla^2E_0\vert\right]\\
&\relphantom{=}+s^{10-8\alpha}\cdot\left[\vert \nabla^2 E\vert\cdot\vert\nabla E_0\vert^3\right].
\end{split}
\end{equation}

Recall that 
\begin{align*}
\frac{\partial^2 \gamma}{\partial q^i \partial q^j}(s,q,p)&=\frac{\partial^2\sigma}{\partial w^i \partial w^j} +\frac{\partial^2 \sigma}{\partial w^i \partial z^k} \left(-\lambda \frac{s^{2-2\alpha}}{2-2\alpha}(\partial_{jk} \phi_0)(q-sp)\right)\\
&\relphantom{=}+\frac{\partial^2 \sigma}{\partial z^k \partial w^l}\left(-\lambda \frac{s^{2-2\alpha}}{2-2\alpha}(\partial_{ik} \phi_0)(q-sp)\right)\\
&\relphantom{=}+\frac{\partial^2 \sigma}{\partial z^k \partial z^l}\left(-\lambda \frac{s^{2-2\alpha}}{2-2\alpha}(\partial_{ik} \phi_0)(q-sp)\right) \left(-\lambda \frac{s^{2-2\alpha}}{2-2\alpha}(\partial_{jl} \phi_0)(q-sp)\right)\\
&\relphantom{=}+\frac{\partial \sigma}{\partial z^k}\left(-\lambda \frac{s^{2-2\alpha}}{2-2\alpha}(\partial_{ijk} \phi_0)(q-sp)\right)
\end{align*}
and so
\begin{align*}
|\nabla_q \gamma|&\apprle |\nabla_w\sigma| + s^{2-2\alpha}|\nabla_z\sigma||\nabla E_0|,\\
|\nabla_q^2 \gamma|&\apprle |\nabla_w^2 \sigma|+s^{2-2\alpha} |\nabla_{w,z}^2\sigma| |\nabla E_0|+s^{4-4\alpha}|\nabla_z^2 \sigma | |\nabla E_0|^2\\
&\relphantom{=}+s^{2-2\alpha} |\nabla_z\sigma| |\nabla^2 E_0|.
\end{align*}
{{Then it follows that 
\begin{equation*}
\begin{split}
\nabla^2{\bf j}(q,s)&=\nabla^2\int_{\mathbb{R}^3} p\gamma^2(s,q,p)dp=\int  2p [(\nabla_q \gamma)^2+ \gamma\nabla_q^2\gamma]dp     \\
\Vert\nabla^2{\bf j}(s)\Vert_{L^\infty_q} 
&\lesssim \norm{\action{p}^{2.1} \nabla_q\gamma}{L^\infty_{q,p}}^2+\norm{\action{p}^5\gamma}{L^\infty_{q,p}}\norm{\nabla_q^2\gamma}{L^\infty_{q,p}}    \\
&\apprle  c_0^2s^{14-14\alpha} \left[\norm{\action{z}^5 \sigma}{\Leb{\infty}_{w,z}}\norm{\nabla^2_{w,z}\sigma}{\Leb{\infty}_{w,z}}+\norm{\action{z}^{2.1}\nabla_{w,z}\sigma}{\Leb{\infty}_{w,z}}^2\right]
\end{split}
\end{equation*}}}
Hence by {{Proposition \ref{prop:continuity-of-E-2}}} and the bootstrap assumptions, we get
\begin{equation*}
\begin{split}
\Vert\nabla^2E(s,q)-\nabla^2E_0(w)\Vert_{L^\infty_{w,z}}&\apprle {{c_0^2 s^{14-14\alpha+\frac{3-\alpha}{2}} =c_0^2 s^{\frac{31-29\alpha}{2}}.}}
\end{split}
\end{equation*}
Therefore, the estimate \eqref{3DerK} follows by \eqref{eq:K-third-1} and \eqref{eq:K-third-2}. 
 This completes the proof of Lemma \ref{lem:particle-density-bootstrap}.
\end{proof}

\subsection{Proof of Theorem \ref{thm:Wave}}\label{subs:wave-theorem}
This subsection is devoted to the proof of Theorem \ref{thm:Wave}. We construct a solution via Picard iteration highly relying on the Hamiltonian structure of the system.

\begin{proof}[Proof of Theorem \ref{thm:Wave}]
Define $\sigma_{(0)}(s,w,z)=\sigma_0(w,z)$ and given $\sigma_{(n)} \in C([0,T];C^1_{w,z})$ satisfying \eqref{eq:wave-bootstrap} with $A_1=A_2=A_3=4c_0$, let $\sigma_{(n+1)} \in C([0,T];C^1_{w,z})$ be the solution of 
\begin{equation}\label{eq:sigma-n-approx}
\begin{aligned}
&\partial_s \sigma_{(n+1)}+\{\sigma_{(n+1)},\mathcal{K}_n\}=0,\quad \sigma_{(n+1)}(0)=\sigma_0, \\
&\mathcal{K}_n:=-\lambda s^{1-2\alpha} (\phi_0(w)-\phi_n(s,q)),\\
& \phi_n(s,q):= \int_{\mathbb{R}^3} V_\alpha(q-y) \gamma_{(n)}^2(s,y,p)\myd{ydp},
\end{aligned}
\end{equation}
{{where $\gamma_{(n)}$}} is defined by $\sigma_{(n)}$ through \eqref{eq:sigma-gamma}. Such a solution exists since \eqref{eq:sigma-n-approx} is {{a linear transport equation.}}  Moreover, by a standard argument, one can show that $\nabla_{z,w}^2 \sigma_{(n)}\in\Leb{\infty}_{s,w,z}$ (see e.g. \cite[Section IV]{BD85}).

By Lemmas \ref{lem:modified-a-priori-electric-field} (or Lemma \ref{lem:a-priori-electric-field} if $\alpha<1$) and \ref{lem:particle-density-bootstrap}, there exists $T(c_0)>0$ such that \eqref{eq:wave-bootstrap} holds for all $s\in [0,T(c_0)]$ with {$A_1=A_2=A_3=c_0$.} Moreover, one can show that \eqref{eq:sigma-c0} and \eqref{eq:sigma-sigma0} hold for $\sigma_{(n)}$ uniformly in $n$ by using the commutation relations given in Subsection \ref{subs:comm-relationship}.

Next, we show that $\sigma_{(n)}$ forms a Cauchy sequence in $\Leb{\infty}_{s,w,z}$. Define  
\[ \delta_{(n)}=\sigma_{(n+1)}-\sigma_{(n)},\quad \delta \mathcal{K}_{(n)}=\mathcal{K}_n-\mathcal{K}_{n-1},\quad \mathfrak{L}_n :=\partial_s + \{\cdot,\mathcal{K}_n\},\quad \delta \mathcal{L}_n=\{\cdot,\delta\mathcal{K}_{(n)}\} \]
so that 
\begin{equation}\label{eq:delta-n}
 \mathfrak{L}_n \delta_{(n)}=\delta \mathfrak{L}_n \sigma_{(n)}.
\end{equation}

Note that 
\begin{equation}\label{eq:Kn-z-w}
\begin{aligned}
\nabla_z \delta \mathcal{K}_{(n)}&=-\lambda s^{2-2\alpha}(E_n(s,q)-E_{n-1}(s,q)),\\
\nabla_w \delta \mathcal{K}_{(n)}&=-\lambda s^{1-2\alpha}(E_n(s,q)-E_{n-1}(s,q))-\lambda^2 \frac{s^{4-4\alpha}}{2-2\alpha} (E_n(s,q)-E_{n-1}(s,q))\cdot \nabla E_0(q).
\end{aligned}
\end{equation}

We claim that if $\alpha<1$, then
\begin{equation}\label{eq:delta-n-sequence}
\norm{\nabla_{w,z}\delta \mathcal{K}_{(n)}(s)}{\Leb{\infty}_{w,z}}\apprle c_0 \left( \norm{\delta_{(n-1)}(s)}{\Leb{2}_{w,z}}+ \norm{\delta_{(n-1)}(s)}{\Leb{\infty}_{w,z}}\right)
\end{equation}
Similarly, if $\alpha>1$, then
\begin{equation}\label{eq:delta-n-sequence-modified}
\norm{\nabla_{w,z}\delta \mathcal{K}_{(n)}(s)}{\Leb{\infty}_{w,z}}\apprle {{c_0s^{13(1-\alpha)}}} \left( \norm{\delta_{(n-1)}(s)}{\Leb{2}_{w,z}}+ \norm{\delta_{(n-1)}(s)}{\Leb{\infty}_{w,z}}\right)
\end{equation}

We only prove the case $\alpha>1$ since the case $\alpha<1$ is simpler than this case. Since $\delta_{(n)}$ satisfies \eqref{eq:delta-n}, using the Hamiltonian structure and the uniform estimate \eqref{eq:sigma-c0}, we get 
\begin{align*}
 \norm{\delta_{(n)}(s)}{\Leb{r}_{w,z}}&\apprle \int_0^s \norm{\nabla_{w,z} \sigma_{(n)}}{\Leb{r}_{w,z}} \norm{\nabla_{w,z} \delta \mathcal{K}_{(n)}}{\Leb{\infty}_{w,z}} d\tau\\
 &\apprle c_0^2 \int_0^s \tau^{13(1-\alpha)}\left(\norm{\delta_{(n-1)}(\tau)}{\Leb{2}_{w,z}}+\norm{\delta_{(n-1)}(\tau)}{\Leb{\infty}_{w,z}}\right)d\tau
\end{align*}
for $r\in \{2,\infty\}$. Since $\alpha<14/13$, if we choose $T$ sufficiently small, then one can easy to check that $\sigma_{(n)}$ form a Cauchy sequence in $ \Leb{\infty}_s\Leb{2}_{w,z}\cap \Leb{\infty}_{s,w,z}$, and hence there exists $\sigma \in \Leb{\infty}_s\Leb{2}_{w,z}\cap \Leb{\infty}_{s,w,z}$ such that $\sigma_{(n)}\rightarrow \sigma$. Moreover, by the uniform convergence, we see that $\sigma \in C([0,T];\Leb{2}_{w,z})$.

On the other hand, by an interpolation inequality \eqref{eq:interpolation-multiplicative}, we have
\begin{align*}
\norm{\nabla_{w,z}\delta_{(n)}}{\Leb{\infty}_{w,z}}&\apprle \norm{\delta_{(n)}}{\Leb{\infty}_{w,z}}^{1/2}\left[\norm{\nabla_{w,z}^2 \sigma_{(n+1)}}{\Leb{\infty}_{w,z}}+\norm{\nabla_{w,z}^2 \sigma_{(n)}}{\Leb{\infty}_{w,z}}\right]^{1/2}.
\end{align*}
Since $\nabla^2_{w,z}\sigma_n$ is uniformly bounded in $\Leb{\infty}_{s,w,z}$, it follows that $\sigma_{(n)}$ forms a Cauchy sequence in $C_s C^1_{w,z}$. Moreover, {{by Fatou's lemma or the conservation}, and \eqref{eq:commutation-relations-VR-wave}}, one can show that $\sigma$ satisfies \eqref{eq:sigma-c0} and \eqref{eq:sigma-sigma0}. Since $\sigma_{(n)}$ converges to $\sigma$ in $C_s C^1_{w,z}$ and $\sigma_{(n+1)}$ satisfies \eqref{eq:sigma-n-approx}, we see that $\sigma$ is a solution of \eqref{EqSigma} with $\sigma(0)=\sigma_0$. {{The proof of uniqueness is straightforward and is therefore omitted}}.

Hence it remains for us to show \eqref{eq:delta-n-sequence} and \eqref{eq:delta-n-sequence-modified}.  First, it follows from \eqref{eq:change-of-variable-sigma-gamma} and \eqref{eq:wave-bootstrap-1} that 
\begin{align*}
&\relphantom{=}|E_n(s,q)-E_{n-1}(s,q)|\\
&\apprle \norm{\delta_{(n-1)}(s)}{\Leb{\infty}_{w,z}} \left(\norm{\action{p}^4\gamma_{(n)}}{\Leb{\infty}_{q,p}}+\norm{\action{p}^4\gamma_{(n-1)}}{\Leb{\infty}_{q,p}} \right)\\
&\relphantom{=}+\norm{\delta_{(n-1)}(s)}{\Leb{2}_{w,z}} \left(\norm{\gamma_{(n)}}{\Leb{2}_{w,z}}+\norm{\gamma_{(n-1)}}{\Leb{\infty}_{w,z}} \right)\\
&\apprle {{c_0^{8}}}s^{4(2-2\alpha)} \norm{\delta_{(n-1)}(s)}{\Leb{\infty}_{w,z}} \left(\norm{\action{z}^4\sigma_{(n)}}{\Leb{\infty}_{w,z}}+\norm{\action{z}^4\sigma_{(n-1)}}{\Leb{\infty}_{w,z}} \right)\\
&\relphantom{=}+\norm{\delta_{(n-1)}(s)}{\Leb{2}_{w,z}} \left(\norm{\sigma_{(n)}}{\Leb{2}_{w,z}}+\norm{\sigma_{(n-1)}}{\Leb{\infty}_{w,z}} \right)\\
&\apprle {{c_0^{9}}} s^{4(2-2\alpha)} \norm{\delta_{(n-1)}(s)}{\Leb{\infty}_{w,z}}+ c_0 \norm{\delta_{(n-1)}(s)}{\Leb{2}_{w,z}}
\end{align*}
and so it follows from  \eqref{eq:Kn-z-w} that 
\begin{align*}
|\nabla_z \delta\mathcal{K}_{(n)}(s,q)|&\apprle  c_0 s^{5(2-2\alpha)}\left(\norm{\delta_{(n-1)}(s)}{\Leb{\infty}_{w,z}} +\norm{\delta_{(n-1)}(s)}{\Leb{2}_{w,z}}\right)
\end{align*}

To estimate $\nabla_w \delta \mathcal{K}_{(n)}$, it follows from \eqref{eq:Kn-z-w} and \eqref{eq:E0-control} that  
\begin{equation}\label{eq:delta-K-n}
\begin{aligned}
|\nabla_w \delta \mathcal{K}_{(n)}(s,q)|&\apprle s^{1-2\alpha} |E_n(s,q)-E_{n-1}(s,q)|+s^{4-4\alpha}|E_n(s,q)-E_{n-1}(s,q)| |\nabla E_0(q)|\\
&\apprle s^{1-2\alpha} |E_n(s,q)-E_{n-1}(s,q)| +c_0 s^{4-4\alpha} |E_n(s,q)-E_{n-1}(s,q)|.
\end{aligned}
\end{equation}
The second part is controlled by 
\[ 
  c_0 s^{12(1-\alpha)}\left(\norm{\delta_{(n-1)}(s)}{\Leb{\infty}_{w,z}} +\norm{\delta_{(n-1)}(s)}{\Leb{2}_{w,z}}\right).
\]

 To estimate the bound for the first term, we note that 
 \begin{equation}\label{eq:difference-density}
  \partial_s (\rho_{(n)}-\rho_{(n-1)})+\Div_q (\delta \boldj_n)=0,\quad \rho_{(n)}-\rho_{(n-1)}=\int_{\mathbb{R}^3} \gamma_{(n)}^2-\gamma_{(n-1)}^2 \myd{p}.
 \end{equation}
Following the argument as in the proof of Proposition \ref{prop:continuity-of-E-2}, it follows from \eqref{eq:difference-density} that  
\begin{align*}
&|E_n(s,q)-E_{n-1}(s,q)|\\
&\apprle s^{2-\alpha} \left( \norm{\delta \boldj_n(s)}{\Leb{\infty}_{q}} + \norm{\delta_{n-1}(s)}{\Leb{2}_{w,z}}\left(\norm{\sigma_{(n)}(s)}{\Leb{2}_{z,w}}+\norm{\sigma_{(n-1)}(s)}{\Leb{2}_{z,w}}\right)\right)
\end{align*}
On the other hand, we have
\begin{align*}
 \norm{\delta \mathbf{j}_n(s)}{\Leb{\infty}_q} &\apprle  \norm{\delta_{(n-1)}(s)}{\Leb{\infty}_{w,z}}\left[\norm{\action{p}^5\gamma_{(n)}(s)}{\Leb{\infty}_{q,p}}+\norm{\action{p}^5\gamma_{(n-1)}(s)}{\Leb{\infty}_{q,p}} \right],\\
 &\apprle s^{10-10\alpha} \norm{\delta_{(n-1)}(s)}{\Leb{\infty}_{w,z}} \left(\norm{\action{z}^5 \sigma_{(n)}(s)}{\Leb{\infty}_{w,z}}+\norm{\action{z}^5 \sigma_{(n-1)}(s)}{\Leb{\infty}_{w,z}} \right).
\end{align*}
Hence we have 
\begin{align*}
|E_n(s,q)-E_{n-1}(s,q)|&\apprle c_0 s^{12-11\alpha} \left( \norm{\delta_{n-1}(s)}{\Leb{2}_{w,z}}+\norm{\delta_{n-1}(s)}{\Leb{\infty}_{w,z}}\right)
\end{align*}
and therefore by \eqref{eq:delta-K-n}, we get 
\[ |\nabla_w \delta \mathcal{K}_{(n)}(s,q)|\apprle c_0 s^{13(1-\alpha)}\left( \norm{\delta_{n-1}(s)}{\Leb{2}_{w,z}}+\norm{\delta_{n-1}(s)}{\Leb{\infty}_{w,z}}\right).\]
This completes the proof of Theorem \ref{thm:Wave}. 
\end{proof}

\subsection{Proof of Theorem \ref{thm:B}}\label{subsec:propagation}

{
In this subsection, we give a proof of Theorem \ref{thm:B}. When $\alpha>1$, we can construct a modified wave operator with an arbitrary size on the final data. In \cite{FOPW23}, they applied the classical result due to Lions and Perthame \cite{LP91} to construct a modified wave operator for the Vlasov-Poisson system. 

A similar result of \cite{LP91} also holds for the Vlasov-Riesz system when $\alpha>1$ by following the argument in Lions and Perthame \cite{LP91}. We note that $\nabla E$ is a weakly singular operator, so the problem is less delicate than the classical Vlasov-Poisson system. In terms of regularity, Lafleche and Saffirio \cite[Proposition A.1]{LS23} proved that the solution becomes classical if $1\leq  \alpha<3/2$ and the initial data has higher moments in $x$ and $v$.

The following lemma can be found in Lafleche \cite[Corollary 5.1]{L19}.
\begin{lemma}
Let $\alpha>1$. Suppose that $E\in \Leb{\infty}_{t,x}$ and let $\mu$ be a solution to the Vlasov-Riesz system with the initial data $\mu_0 \in \Leb{2}_{x,v}\cap \Leb{\infty}_{x,v}$, and $\mu_0^2 |v|^n \in \Leb{2}_{x,v}\cap \Leb{\infty}_{x,v}$ with $n>3$. Then $\rho \in \Leb{\infty}_{t}([0,T];L^\infty_x)$. Moreover, we have 
\begin{equation}
\norm{\rho}{\Leb{\infty}_{t,x}}\apprle \left(\norm{\mu^2_0|v|^n}{\Leb{\infty}_{x,v}}^{1/n}+\norm{E}{\Leb{\infty}_{t,x}}\norm{\mu^2_0}{\Leb{\infty}_{x,v}}^{1/n}T \right)^n +\norm{\mu_0}{\Leb{\infty}_{x,v}}^2.
\end{equation}
\end{lemma}
%

As an application, since $\alpha>1$, $\nabla_x E$ is written as a weakly singular integral operator of $\rho$. So we can show that 
\begin{equation}\label{eq:nabla-E-control-easy}
\norm{\nabla E}{\Leb{\infty}_{t,x}}\apprle \norm{\rho}{\Leb{\infty}_t\Leb{1}_x}+\norm{\rho}{\Leb{\infty}_{t,x}}.
\end{equation}
This is sufficient for showing $C^1$-regularity of solutions to the Vlasov-Riesz system by following the argument in Lafleche and Saffirio \cite[Lemma A.2]{LS23}.

In summary, we have obtained the following. Since the argument is essentially the same, we do not reproduce the proof here.

\begin{theorem}\label{thm:propagation}
There exists $\delta>0$ such that for $\alpha \in (1,1+\delta)$, if $\mu_0 \in \Leb{2}_{x,v}\cap \Leb{\infty}_{x,v}$ satisfies 
\[ \iint_{\mathbb{R}^3_x\times\mathbb{R}^3_v} |v|^m \mu_0^2(x,v)dxdv<\infty \]
for all $m<m_0$ and $m_0>3$, then there exists a solution $\mu\in C(\mathbb{R};\Leb{p}(\mathbb{R}^3_x\times\mathbb{R}^3_v))\cap \Leb{\infty}(\mathbb{R}_t\times\mathbb{R}^3_x\times\mathbb{R}^3_v)$ for all $2\leq p<\infty$ of \eqref{eq:VR-mu} satisfying
\begin{equation}
\sup_{t\in[0,T]} \iint_{\mathbb{R}^3_x\times\mathbb{R}^3_v} |v|^m \mu^2(t,x,v)dxdv<\infty 
\end{equation}
for all $T<\infty$ and $m<m_0$. Moreover, if  $\mu_0\in \Leb{2}_{x,v}\cap\Leb{\infty}_{x,v}$ satisfies 
\[ \norm{\action{v}^{m} \mu_0}{\Leb{2}_{x,v}\cap \Leb{\infty}_{x,v}}+\norm{\nabla_{x,v}\mu_0}{\Leb{2}_{x,v}\cap\Leb{\infty}_{x,v}}<\infty\]
for some $m>3$ in addition, then $E\in \Leb{\infty}_{t,x}$ and $\nabla_{x,v}\mu \in\Leb{\infty}_{t,x,v}$.
\end{theorem}

\begin{proof}[Proof of Theorem \ref{thm:B}]
(i) Note that the assumption on $\mu_\infty$ and $E_\infty$ in \ref{eq:assumption-wave} implies \eqref{eq:E0-control} and \eqref{eq:sigma0-control}. By Theorem \ref{thm:Wave}, there exists a unique $\sigma\in C([0,T^*];\Leb{2}_{w,z})$  of \eqref{EqSigma}. Moreover, $\sigma$ satisfies \eqref{eq:sigma-c0} and \eqref{eq:sigma-sigma0}. Suppose that $\alpha>1$. If we use \eqref{eq:change-of-variable-sigma-gamma} and the pseudo-conformal transform again,  then we can apply Theorem \ref{thm:propagation} together with the time-reversal symmetry to construct a modified wave operator, which proves (i)-(a).

In addition, if we assume the smallness on the final data, then there exists $\varepsilon_0>0$ such that there exists a unique $\sigma\in C([0,1];\Leb{2}_{w,z})$  of \eqref{EqSigma}. Moreover, $\sigma$ satisfies \eqref{eq:sigma-c0} and \eqref{eq:sigma-sigma0}. From this, if we use \eqref{eq:change-of-variable-sigma-gamma} and the pseudo-conformal transform again, then if we choose $\varepsilon_0>0$ sufficiently small, then it follows from Theorem \ref{thm:VR-LWP} that the problem admits a local solution from $t=1$ to $t=0$ by using the time-reversal symmetry. This completes the proof of (i)-(b).

(ii)  We first apply Theorem \ref{thm:B} (i) and then Theorem \ref{thm:A}.
\end{proof}

\begin{remark}\label{rem:wave-operator-large}
If $\alpha>1$, we applied a similar result to \cite{LP91} to construct the modified wave operator. It seems that it is delicate to obtain a corresponding result when $\alpha<1$ since $\nabla E$ is a hypersingular integral operator of $\rho$ which is difficult to bound in $\Leb{\infty}_x$.
\end{remark}
}

\section*{Declarations}
\subsection*{Conflict of interest}
 The authors have no relevant financial or non-financial interests to disclose.

\appendix

\section{Local well-posedness of strong solutions to Vlasov-Riesz system}\label{app:A}

In this section, we prove the local well-posedness of the Vlasov-Riesz system of order $\alpha$ for any $1/2<\alpha<3/2$. We remark that the proof of the local well-posedness of the system in this regime seems to be classical, but we cannot find an appropriate reference to cite it that we need to use. Hence, for the sake of completeness, we provide proof by using Picard iteration.

We recall the Vlasov-Riesz system of order $\alpha$:
\begin{equation}\label{eq:VR-Hamiltonian}
\partial_t \mu + v \cdot \nabla_x \mu - \nabla_x \phi \cdot \nabla_v \mu =0,
\end{equation}
where 
\[ \phi(t,x)=c_\alpha \int_{\mathbb{R}^3_{y,v}} \frac{\mu^2(t,y,v)}{|x-y|^{3-2\alpha}} dydv. \] 

If we define the Poisson bracket by 
\[ \{f,g\}=\nabla_x f \cdot \nabla_v g -\nabla_x g \cdot \nabla_v f,\]
then one can rewrite \eqref{eq:VR-Hamiltonian} into
\begin{equation*}\label{eq:VR-Hamiltonian-Strong}
 \partial_t \mu + \{\mu,\mathcal{H}\}=0,\quad \mathcal{H}=\frac{|v|^2}{2}- c_\alpha \iint_{\mathbb{R}^3_{y,v}}  \mu^2(t,y,v)\frac{dydv}{|x-y|^{3-2\alpha}}.
\end{equation*}

For $0<T\leq \infty$, we say that $\mu \in X_T$ if $\mu \in \Leb{2}((0,T)\times \mathbb{R}^3_x\times \mathbb{R}^3_v)$ satisfies $\partial_t \mu + v\cdot \nabla_x \mu$, $\nabla_x\mu,\nabla_v \mu \in \Leb{2}((0,T)\times \mathbb{R}^3_x\times \mathbb{R}^3_v)$.

We say that $\mu$ is a strong solution to the Vlasov-Riesz system if $\mu\in X_T$ and  $\mu \in C([0,T];\Leb{2}_{x,v})$ satisfies  \eqref{eq:VR-Hamiltonian} pointwise a.e. in $(0,T)\times \mathbb{R}^3_x\times \mathbb{R}^3_v$.

\begin{theorem}\label{thm:VR-LWP}
Let $1/2<\alpha<3/2$ and let $\mu_0$ satisfy 
\[ B:=\norm{\mu_0}{\Leb{2}_{x,v}}+\norm{\action{x,v}^4\mu_0}{\Leb{\infty}_{x,v}}+\norm{\nabla_{x,v}\mu_0}{\Leb{2}_{x,v}\cap \Leb{\infty}_{x,v}}<\infty.\]
Then there exist $T=T(\alpha,B)>0$ and a constant $C>0$ depending only on $\alpha$ such that the problem \eqref{eq:VR-Hamiltonian} has a unique strong solution $\mu \in X_T$  satisfying
\[ \sup_{t\in [0,T]}\norm{\mu(t)}{\Leb{2}_{x,v}}= \norm{\mu_0}{\Leb{2}_{x,v}},\quad \sup_{t\in [0,T]}\norm{\action{x,v}^4\mu(t)}{\Leb{\infty}_{x,v}}+\norm{\nabla_{x,v}\mu(t)}{\Leb{2}_{x,v} \cap \Leb{\infty}_{x,v}}\leq CB. \]
Moreover, if $\mu_0 \in C^1_{x,v}$ in addition, then the strong solution becomes classical.
\end{theorem} 
{
\begin{remark}
By \eqref{eq:T-choice-1} and \eqref{eq:T-choice-2}, one can see that if $B$ is sufficiently small, then $T>1$. 
\end{remark}
}

\begin{proof}
Define the following iteration scheme:
\[ \phi_0(t,x)=0,\quad \mu_0(t,x,v)=\mu_0(x,v) \]
\begin{equation}\label{eq:Picard-VR}
 \partial_t \mu_{n+1}+\{\mu_{n+1},\mathcal{H}_n\}=0,\quad \mu_{n+1}(t=0)=\mu_0,\quad \mathcal{H}_n:=\frac{|v|^2}{2}- \phi_n(t,x),
\end{equation}
where 
\[ \phi_n(t,x)=c_\alpha \iint_{\mathbb{R}^3_{y,v}} \frac{\mu_n^2(t,y,v)}{|x-y|^{3-2\alpha}}dydv.\]

It is easy to see that
\begin{equation}\label{eq:gradient-phi-n}
   \norm{D\phi_n}{\Leb{\infty}_x}\apprle \norm{\mu_n}{\Leb{2}_{x,v}}^2 + \norm{\action{v}^2 \mu_n}{\Leb{\infty}_{x,v}}^2.
\end{equation}
\begin{equation}\label{eq:hessian-phi-n}
\norm{D^2 \phi_n}{\Leb{\infty}_x}\apprle  \norm{\action{v}^4\mu_n}{\Leb{\infty}_{x,v}} \norm{\nabla_x \mu_n}{\Leb{\infty}_{x,v}}+ \norm{\mu_n}{\Leb{2}_{x,v}}\norm{\nabla_x \mu_n}{\Leb{2}_{x,v}}.
\end{equation} 

If we write 
\[ V_n=(v,-\nabla_x \mathcal{H}_n),\]
then $\Div_{x,v} V_n=0$ and 
\begin{equation}
 \partial_t \mu_{n+1}+\Div(\mu_{n+1}V_n)=0,\quad \mu_{n+1}(t=0)=\mu_0,
\end{equation}
{which is a linear transport equation. Hence, by the Cauchy-Lipschitz theorem (see e.g. \cite[Theorem 3.1]{MB01}), the problem \eqref{eq:Picard-VR} admits a unique solution $\mu_{n+1} \in C^1([0,\infty); \Sob{1}{\infty}_{x,v})$. Moreover, since $\Div_{x,v} V_n=0$, it follows from Liouville's theorem that 
\[ \norm{\mu_{n+1}(t)}{\Leb{p}_{x,v}}=\norm{\mu_0}{\Leb{p}_{x,v}} \]
for all $t\in [0,\infty)$ and $1\leq p\leq \infty$.}

For $\omega=\omega(v)$ or {$\omega=\omega(x)$},  a direct computation gives 
\begin{equation}\label{eq:weight-v}
\partial_t (\omega \mu_{n+1})+\{\omega \mu_{n+1},\mathcal{H}_n\}=\mu_{n+1}\{\omega,\mathcal{H}_n\}.
\end{equation}
If we choose $\omega(v)=1$ in \eqref{eq:weight-v}, then we have
\[
   \norm{\mu_{n+1}}{\Leb{2}_{x,v}} = \norm{\mu_0}{\Leb{2}_{x,v}}.
\]
Also, by formally differentiate \eqref{eq:Picard-VR} in $x$ and $v$, we have
\begin{equation}\label{eq:D-derivative}
\begin{aligned}
\partial_t (D_{x_i} \mu_{n+1})+\{D_{x_i}\mu_{n+1},\mathcal{H}_n\} &= -\nabla_x(\partial_{x_i}\phi)\cdot \nabla_v \mu_{n+1},\\
\partial_t (D_{v_i} \mu_{n+1})+\{D_{v_i}\mu_{n+1},\mathcal{H}_n\} &= - D_{x_i} \mu_{n+1}.
\end{aligned}
\end{equation}
 for $i=1,2,3$. Hence by \eqref{eq:weight-v}, \eqref{eq:D-derivative}, and Liouville's theorem, for $r\in\{2,\infty\}$, we get 
 \begin{equation}\label{eq:A-estimates}
\begin{aligned}
\norm{\action{v}^a \mu_{n+1}(t)}{\Leb{r}_{x,v}}&\leq \norm{\action{v}^a \mu_{n+1}(0)}{\Leb{r}_{x,v}}+\int_0^t\norm{\nabla_x \phi}{\Leb{\infty}_x} \norm{\action{v}^a \mu_{n+1}}{\Leb{r}_{x,v}}ds,\\
\norm{\action{x}^a \mu_{n+1}(t)}{\Leb{r}_{x,v}}&\leq \norm{\action{x}^a \mu_{n+1}(0)}{\Leb{r}_{x,v}}+\int_0^t\norm{\nabla_x \phi}{\Leb{\infty}_x} \norm{\action{x}^a \mu_{n+1}}{\Leb{r}_{x,v}}ds,\\
\norm{D_{x_j}\mu_{n+1}(t)}{\Leb{r}_{x,v}}&\leq \norm{\nabla_x\mu_0}{\Leb{r}_{x,v}}+\int_0^t \norm{\nabla_x^2\phi_n}{\Leb{\infty}}\norm{\nabla_v \mu_{n+1}}{\Leb{r}_{x,v}}ds,\\
\norm{D_{v_j}\mu_{n+1}(t)}{\Leb{r}_{x,v}}&\leq \norm{\nabla_v\mu_0}{\Leb{r}_{x,v}}+\int_0^t \norm{\nabla_x \mu_{n+1}}{\Leb{r}_{x,v}}ds.
\end{aligned}
\end{equation}
The above calculation can be justified by using a method of finite difference quotient. 

{
If we define 
\begin{align*}
 A_n(t)&=\norm{\nabla_{x,v} \mu_{n}(t)}{\Leb{2}_{x,v}\cap\Leb{\infty}_{x,v}}+\norm{\action{x}^4 \mu_{n+1}(t)}{\Leb{\infty}_{x,v}}+\norm{\action{v}^4 \mu_{n+1}(t)}{\Leb{\infty}_{x,v}},\\
A_0 &= \norm{\nabla_{x,v}\mu_0}{\Leb{2}_{x,v}\cap\Leb{\infty}_{x,v}}+\norm{\action{x}^4\mu_0}{\Leb{\infty}_{x,v}}+\norm{\action{v}^4\mu_0}{\Leb{\infty}_{x,v}},
\end{align*}
then it follows from \eqref{eq:A-estimates} that 
\[ A_{n+1}(t)\apprle A_0 + \int_0^t (1+\norm{\nabla_x \phi_n}{\Leb{\infty}_x}+\norm{\nabla^2_x \phi_n}{\Leb{\infty}_x})A_{n+1}(s)\myd{s}.\]
Hence by Gronwall's inequality, we get 
\begin{align*}
A_{n+1}(t)&\apprle A_0 \exp\left(\int_0^t (1+\norm{\nabla_x \phi_n}{\Leb{\infty}_x}+\norm{\nabla_x^2\phi_n}{\Leb{\infty}_x})ds\right).
\end{align*}

If we define 
\[ B_n(t)=\norm{\mu_n(t)}{\Leb{2}_{x,v}}+A_n(t), \]
then there exist $C>1$ and $T>0$ such that 
\[
B_n(t)\leq CB 
\]
for all $t\in [0,T]$ and $n$. The case $n=0$ is clear. Note that 
\begin{align*}
&B_{n+1}(t)^2\\
&=(\norm{\mu_0}{\Leb{2}_{x,v}}+A_{n+1}(t))^2\\
&\leq 2\norm{\mu_0}{\Leb{2}_{x,v}}^2+2 C_0^2 A_0^2 \exp \left(2\int_0^t (1+\norm{\nabla_x \phi_n}{\Leb{\infty}_x}+\norm{\nabla^2 \phi_n}{\Leb{\infty}_x})ds\right)\\
&\leq 2\norm{\mu_0}{\Leb{2}_{x,v}}^2+2C_0^2 A_0^2 \exp \left(C_1 \int_0^t (1+B_n^2(s))ds\right)\\
&\leq 2C_0^2B^2 \left(1+ \exp\left(C_1 \int_0^t (1+B_n^2(s))ds \right)\right)\\
&\leq 4C_0^2 B^2 \exp\left(C_1 T (1+C^2 B^2)\right).
\end{align*}
Now choose $C=4C_0$ and $T>0$ so that $\exp\left(C_1 T (1+16C_0^2 B^2)\right)\leq 4$. Here 
\begin{equation}\label{eq:T-choice-1}
 4C_1 T (1+16C_0^2 B^2)\leq \ln 4.
\end{equation}
Then we get 
\[ B_{n+1}(t)\leq C B \]
for all $t\in [0,T]$.}

If we write $\delta_n:=\mu_{n+1}-\mu_n$, then a direct computation gives
\begin{equation}\label{eq:Hamiltonian-difference} 
\begin{aligned}
0&=\partial_t \delta_n + \left\{ \delta_n, \frac{\mathcal{H}_n+\mathcal{H}_{n-1}}{2}\right\}- \nabla_x (\phi_n -\phi_{n-1})\cdot \nabla_v \left(\frac{\mu_{n+1}+\mu_n}{2}\right).
\end{aligned}
\end{equation}

Similar to the calculation in \eqref{eq:gradient-phi-n}, we have 
\begin{align*}
\norm{\nabla_x [\phi_n -\phi_{n-1}]}{\Leb{\infty}_x}&\apprle \norm{\delta_{n-1}}{\Leb{\infty}_{x,v}} \left(\norm{\action{v}^4\mu_n}{\Leb{\infty}_{x,v}}+\norm{\action{v}^4\mu_{n-1}}{\Leb{\infty}_{x,v}}\right) \\
&\relphantom{=}+\norm{\delta_{n-1}}{\Leb{2}_{x,v}} \left(\norm{\mu_n}{\Leb{2}_{x,v}}+\norm{\mu_{n-1}}{\Leb{2}_{x,v}}\right)
\end{align*}
Hence by multiplying \eqref{eq:Hamiltonian-difference} with $\delta_n(t)$, we get 
\begin{align*}
\frac{1}{2}\frac{d}{dt}\norm{\delta_n(t)}{\Leb{2}_{x,v}}^2&\apprle \norm{\nabla_x[\phi_n-\phi_{n-1}]}{\Leb{\infty}_{x,v}}\left(\norm{\mu_{n+1}}{\Leb{2}_{x,v}}+\norm{\mu_{n}}{\Leb{2}_{x,v}}\right)\norm{\delta_n(t)}{\Leb{2}_{x,v}}\\
&\apprle B^2 \norm{\delta_{n-1}(t)}{\Leb{2}_{x,v}\cap \Leb{\infty}_{x,v}} \norm{\delta_n(t)}{\Leb{2}_{x,v}}
\end{align*}

By Liouville's theorem induced by the { new Hamiltonian $(\mathcal{H}_n+\mathcal{H}_{n-1})/2$,} we have 
\begin{align*}
\norm{\delta_n(t)}{\Leb{\infty}_{x,v}}&\leq \frac{1}{2}\int_0^t \norm{\nabla_x (\phi_n-\phi_{n-1})}{\Leb{\infty}_{x,v}}(\norm{\nabla_v \mu_{n+1}}{\Leb{\infty}_{x,v}}+\norm{\nabla_v \mu_n}{\Leb{\infty}_{x,v}})ds\\
&\apprle B^2\int_0^t \norm{\delta_{n-1}(s)}{\Leb{2}_{x,v}\cap\Leb{\infty}_{x,v}}ds
\end{align*}
Hence by Gronwall's inequality, we get 
\begin{align*}
\norm{\delta_n(t)}{\Leb{2}_{x,v}\cap \Leb{\infty}_{x,v}}\leq CB^2 \int_0^t \norm{\delta_{n-1}(s)}{\Leb{2}_{x,v}\cap \Leb{\infty}_{x,v}}ds
\end{align*}
for some constant $C>0$. Then by iteration, we get 
\[ \norm{\delta_n(t)}{\Leb{2}_{x,v}\cap \Leb{\infty}_{x,v}}\leq (CB^2T)^n\]
for all $t\in [0,T]$ and for all $n$. This implies that for any $n,m\geq 1$, we have 
\begin{align*}
\norm{\mu_{m}(t)-\mu_{n}(t)}{\Leb{2}_{x,v}\cap \Leb{\infty}_{x,v}}&\leq \sum_{k=n+1}^{m-1} (CB^2T)^k
\end{align*}
for all $t\in [0,T]$.  By choosing $T>0$ sufficiently small so that 
\begin{equation}\label{eq:T-choice-2}
CB^2T\leq 1/2, 
\end{equation}
we see that the sequence $\{\mu_n\}$ is Cauchy in $\Leb{\infty}_t\Leb{2}_{x,v}\cap \Leb{\infty}_{t,x,v}$.  Then by a standard compactness argument and duality argument, we can show the existence of a strong solution satisfying $\mu \in C([0,T];\Leb{2}_{x,v})$. The uniqueness part is also easy to show. Finally, the regularity result follows from a method of characteristics induced by the electric field generated by $\mu$ and the uniqueness of the strong solution.  This completes the proof of Theorem \ref{thm:VR-LWP}.  
\end{proof}

\section{Interpolation inequalities}\label{app:interpolation}

We prove the following interpolation inequalities.
\begin{proposition}\label{prop:interpolation}
Let $s\geq 0$ and $f:\mathbb{R}^3_{x}\times \mathbb{R}^3_v \rightarrow \mathbb{R}$. 
\begin{enumerate}[label=\textnormal{(\roman*)}]
\item for $s\geq 0$, we have 
\begin{equation}
\norm{\action{x}^s \nabla_v f}{\Leb{\infty}_{x,v}}\apprle \norm{\action{x}^{2s} f}{\Leb{\infty}_{x,v}}+\norm{\nabla_v^2 f}{\Leb{\infty}_{x,v}}\label{eq:interpolation-1}.
\end{equation}
In particular, when $s=0$, we have 
\begin{equation}
 \norm{D f}{\Leb{\infty}_{x,v}}\apprle \norm{f}{\Leb{\infty}_{x,v}}^{1/2}\norm{D^2f}{\Leb{\infty}_{x,v}}^{1/2},\quad D \in \{\nabla_x ,\nabla_v\}\label{eq:interpolation-multiplicative}
\end{equation}
\item if $3/2<s<4$, then
\begin{equation}
\norm{\action{x}^s \nabla_x f}{\Leb{\infty}_{x,v}}\apprle \norm{\action{x}^{2s} f}{\Leb{\infty}_{x,v}}+\norm{\nabla_x^2 f}{\Leb{\infty}_{x,v}}.\label{eq:interpolation-2}
\end{equation} 
\end{enumerate}
\end{proposition}

We believe that inequalities should be well-known, but we cannot find an exact reference that contains these inequalities. Hence we provide their proof for completeness by using Littlewood-Paley projections.

 Let $\varphi \in C_c^\infty (B_2)$ be a symmetric radial bump function satisfying $0\leq \varphi \leq1$ and $\varphi = 1$ in $B_1$. Define $\psi(\xi)=\varphi(\xi)-\varphi(2\xi)$. For a dyadic number $N$, define $\varphi_N(\xi)=\varphi(\xi/N)$ and $\psi_N(\xi)=\psi(\xi/N)$. Then we have
 \[ \varphi(\xi)+\sum_{N \geq 1} \psi_N(\xi)=1,\quad \xi \in \mathbb{R}^3.\]
We define the Littlewood-Paley projection operators as Fourier multiplier operators:
\begin{align*}
\widehat{P_{\leq N} f}(\xi)&=\varphi_N(\xi)\hat{f}(\xi),\\
\widehat{P_{> N} f}(\xi)&=(1-\varphi_N(\xi))\hat{f}(\xi),\\
\widehat{P_{N} f}(\xi)&=\psi_N(\xi)\hat{f}(\xi).
\end{align*}
Note that 
\begin{equation}\label{eq:projection-convolution}
P_{\leq N} f(x)=(\check{\varphi}_N*f)(x),\quad \text{where } \check{\varphi}_N(x)=N^3 \check{\varphi}(Nx).
\end{equation}

We recall several properties on Littlewood-Paley projectors (see e.g. \cite[Appendix A]{T06}).
\begin{proposition}\label{prop:LP}
Let $1\leq p\leq q\leq \infty$, $k\in \mathbb{N}$, and let $N$ be a dyadic integer. Then
\begin{enumerate}[label=\textnormal{(\roman*)}]
\item $\norm{P_{\leq N} f}{\Leb{p}}+\norm{P_N f}{\Leb{p}}\apprle \norm{f}{\Leb{p}}$.
\item $N^k \norm{P_N f}{\Leb{p}}\apprle \norm{\nabla^k f}{\Leb{p}}$
\item $\norm{P_N (\nabla^k f)}{\Leb{q}}\apprle N^{k(1/p-1/q)} \norm{f}{\Leb{p}}$.
\end{enumerate} 
\end{proposition}

\begin{proof}[Proof of Proposition \ref{prop:interpolation}]
We first show \eqref{eq:interpolation-1}. Let us take the Littlewood-Paley projection in $v$ variable. Then 
\begin{align*}
\nabla_v f&=\nabla_v P_{\leq N} f+ \nabla_v (f-P_{\leq N} f)
\end{align*}
and so 
\begin{align*}
\norm{\action{x}^s\nabla_v f}{\Leb{\infty}_v}&\apprle \action{x}^s N \norm{f}{\Leb{\infty}_v}+ \action{x}^s \norm{\nabla_v (f-P_{\leq N} f)}{\Leb{\infty}_v}.
\end{align*}

Since $\varphi_N(0)=1$, it follows that 
\begin{align*}
\nabla_v (f-P_{\leq N} f)(x,v)&=\nabla_v \int_{\mathbb{R}^3_v} \check{\varphi}_N(w)[f(x,v)-f(x,v-w)]dw.
\end{align*}
By the mean value theorem and a change of variable, we get
\begin{align*}
|\nabla_v (f-P_{\leq N} f)(x,v)|&\apprle \norm{\nabla_v^2 f(x,\cdot)}{\Leb{\infty}_v} \int_{\mathbb{R}^3_v} |w| | \check{\varphi}_N(w)| dw\apprle N^{-1} \norm{\nabla_v^2 f}{\Leb{\infty}_v}.
\end{align*}
Hence it follows that 
\begin{align*}
\left|\action{x}^s \nabla_v f(x,v)\right|&\apprle \action{x}^s N \norm{f}{\Leb{\infty}_{x,v}}+ \action{x}^sN^{-1}\norm{\nabla_v^2 f(x,\cdot)}{\Leb{\infty}_v}.
\end{align*}
By choosing $N\approx \action{x}^s$, we get 
\[ \norm{\action{x}^s f}{\Leb{\infty}_{x,v}}\apprle \norm{\action{x}^{2s} f}{\Leb{\infty}_{x,v}}+\norm{\nabla_v^2 f}{\Leb{\infty}_{x,v}}.\]
This completes the proof of \eqref{eq:interpolation-1}.

If $s=0$, then the desired inequality \eqref{eq:interpolation-multiplicative} follows by choosing $N\approx \norm{\nabla^2_v f}{\Leb{\infty}_{x,v}}/ \norm{f}{\Leb{\infty}_{x,v}}$.

To show \eqref{eq:interpolation-2}, due to \eqref{eq:interpolation-multiplicative}, it suffices to show that 
\begin{equation}\label{eq:As-estimate}
   \sup_{A\geq 2^{10}} A^s \norm{\chi(A^{-1} x)\nabla_x f}{\Leb{\infty}_x} \apprle \norm{\action{x}^{2s}f}{\Leb{\infty}_x} + \norm{\nabla_x^2 f}{\Leb{\infty}_x}
\end{equation}
since 
\begin{align*}
\norm{\action{x}^s\nabla_x f}{\Leb{\infty}_x}&\apprle \norm{\nabla_x f}{\Leb{\infty}_x} + \norm{|x|^{s} \nabla_x f}{\Leb{\infty}_x(\{|x|\geq 2^{10}\})}\\
&\apprle  \norm{f}{\Leb{\infty}_x}+\norm{\nabla^2_x f}{\Leb{\infty}_x} +\sup_{A\geq 2^{10}} \norm{|x|^{s} \nabla_x f}{\Leb{\infty}_x(\{|x| \in [A,2A]\})}\\
&\apprle \norm{\action{x}^{2s} f}{\Leb{\infty}_x}+\norm{\nabla_x^2 f}{\Leb{\infty}_x}.
\end{align*}

Decompose
\[ f = P_{\leq N_0} f + P_{>N_0} f,\]
where $N_0$ will be chosen later. By \eqref{eq:projection-convolution}, we have
\begin{equation}
|\nabla P_{\leq N_0} f(x)|\apprle N_0^{4} \norm{f}{\Leb{1}}\apprle N_0^4 \norm{\action{x}^{2s} f}{\Leb{\infty}_x},
\end{equation}
which implies 
\begin{equation}\label{eq:low-frequency-1}
A^s \norm{\chi(A^{-1}x)\nabla_x P_{\leq N_0} f}{\Leb{\infty}_x} \apprle A^s N_0^4 \norm{\action{x}^{2s} f}{\Leb{\infty}_x}.
\end{equation}
For $A> 2^{10}$, we can choose a dyadic number $N_0$ so that $A^s N_0^4\leq 1$ and $A^{{s}/{3}}N_0 >1$. So now we estimate 
\begin{align*}
 \sum_{N\geq N_0} A^s \norm{\chi(A^{-1}x)\nabla_x P_{N} f}{\Leb{\infty}_x}&\leq  \sum_{A^{-s/3}\leq N \leq A^s } A^s \norm{\chi(A^{-1}x)\nabla_x P_{N} f}{\Leb{\infty}_x}\\
 &\relphantom{=}+\sum_{N\geq A^s} A^s \norm{\chi(A^{-1}x)\nabla_x P_{N} f}{\Leb{\infty}_x}.
\end{align*}

To estimate the contribution due to a higher frequency of $f$, we rewrite
\begin{align*}
\nabla P_N f&= N Q_N f,
\end{align*}
where $\widehat{Q_N f}(\xi) = i\frac{\xi}{N} \psi_N(\xi) \hat{f}(\xi)$. By a change of variable, we can rewrite
\begin{equation}\label{eq:convolution-QN}
Q_N f(x)=\int_{\mathbb{R}^3} N^3 m(N(x-y))f(y)dy,
\end{equation}
where 
\[  m(x)= \int_{\mathbb{R}^3} i \xi \psi(\xi)e^{ix\cdot \xi} d\xi.\]
Then by \eqref{eq:convolution-QN}, we have
\[ \norm{Q_Nf}{\Leb{p}}\apprle \norm{f}{\Leb{p}},\quad 1\leq p\leq \infty.\]
Moreover, the operator $Q_N$ shares similar properties as in Littlewood-Paley projectors.

So for $N\geq A^s$, it follows from Proposition \ref{prop:LP} that 
\begin{equation*}
\begin{aligned}
  A^s \norm{\chi(A^{-1} x)\nabla P_N f(x)}{\Leb{\infty}_x}&\apprle A^s N\norm{\chi(A^{-1}x) Q_Nf(x)}{\Leb{\infty}_x}\apprle A^s N^{-1} \norm{\nabla^2 f}{\Leb{\infty}_x}.
  \end{aligned}
 \end{equation*}
 Hence we conclude that
  \begin{equation}\label{eq:high-frequency-1} 
 \sum_{N\geq A^s} A^s N^{-1} \norm{\nabla^2 f}{\Leb{\infty}_x}\apprle \norm{\nabla^2 f}{\Leb{\infty}_x}. 
 \end{equation}
 
 To estimate the intermediate pieces $A^{-s/3}\leq N \leq A^s$, for sufficiently small $c \ll 1$ and $\vartheta$ satisfying $\sum_{p \in \mathbb{Z}^3} \vartheta(x-p)=1$, we rewrite
 \begin{align*}
Q_Nf(x) &=\sum_{p \in \mathbb{Z}^3} \vartheta(c^{-1} A^{-1}x-p)\int_{\mathbb{R}^3} N^3 m(N(x-y)) f(y)dy=\sum_{p,q\in \mathbb{Z}^3} T_{p,q}f(x),
 \end{align*}
 where
\begin{align*}
T_{p,q}^Nf(x)&=\int_{\mathbb{R}^3} K_{p,q}^N(x,y) f(y)dy,\\
  K_{p,q}^N(x,y)&=\vartheta(c^{-1} A^{-1}x-p)\cdot N^3 m(N(x-y)) \cdot \vartheta (c^{-1} A^{-1} y-q).
 \end{align*}

Note that 
 \begin{equation}\label{eq:TN-pq}
 \norm{T_{p,q}^N}{\Leb{\infty}_x\rightarrow\Leb{\infty}_x}\apprle (cNA)^{-10} |p-q|^{-10}
 \end{equation}
 whenever $|p-q|\geq 100$. Indeed, since $\supp \vartheta \subset [-1/2,1/2]^3$ and $|p-q|\geq 100$, we have $|x-y|\geq \frac{99}{100} cA|p-q|$. This implies
 \[
   |K_{p,q}^N(x,y)|\apprle \frac{N^3}{(1+N|x-y|)^{100}}\apprle  \frac{N^3}{(1+N|x-y|)^{10}} \cdot \frac{1}{(cNA|p-q|)^{10}}
\]
 
 Hence it follows that 
 \begin{align*}
 A^s N \norm{\chi(A^{-1} x)Q_N f}{\Leb{\infty}_x}&\apprle A^s N \sum_{p,q\in \mathbb{Z}^3} \norm{\chi(A^{-1} x)T_{p,q}^Nf}{\Leb{\infty}_x}\\
 &\apprle A^s N \sum_{c^{-1}/10 \leq |p|\leq 10 c^{-1}} \sum_{|p-q|\leq 100}  \norm{\chi(A^{-1} x) T_{p,q}^Nf}{\Leb{\infty}_x}\\
 &\relphantom{=}+ A^s N \sum_{c^{-1}/10 \leq |p|\leq 10 c^{-1}} \sum_{|p-q|\geq 100}  \norm{\chi(A^{-1} x)T_{p,q}^Nf}{\Leb{\infty}_x}.
 \end{align*}
 
 For the first part, we have
 \begin{align*}
  A^s N \sum_{c^{-1}/10 \leq |p|\leq 10 c^{-1}} \sum_{|p-q|\leq 100}  \norm{\chi(A^{-1} x)T_{p,q}^Nf}{\Leb{\infty}_x}\apprle A^s N \norm{1_{\{|x|\geq 10^{-5} A\}} f}{\Leb{\infty}_x}.
 \end{align*}
since $\supp \vartheta \subset [-1/2,1/2]^3$, $|p-q|\leq 100$, and $\supp \chi \subset B_4\setminus B_{1/4}$.  For the second part, it follows from \eqref{eq:TN-pq} that
 \begin{align*}
 A^s N \sum_{c^{-1}/10 \leq |p|\leq 10 c^{-1}} \sum_{|p-q|\geq 100}  \norm{\chi(A^{-1} x)T_{p,q}^Nf}{\Leb{\infty}_x}&\apprle A^s N \sum_{q\in \mathbb{Z}^3} \action{q}^{-10} (cNA)^{-10} \norm{f}{\Leb{\infty}_x}\\
&\apprle A^{s-4}  \norm{f}{\Leb{\infty}_x}.
 \end{align*}
Hence 
 \begin{align*}
 \sum_{A^{-s/3}\leq N \leq A^s} A^s N \norm{\chi(A^{-1} x)Q_N f}{\Leb{\infty}_x}&\apprle  \sum_{A^{-s/3}\leq N \leq A^s} A^s N \norm{1_{\{|x|\apprge A\}} f}{\Leb{\infty}_x} + A^{s-4} \norm{f}{\Leb{\infty}_x} \sum_{A^{-s/3} \leq N \leq A^s} 1\\
 &\apprle A^{2s} \norm{1_{\{|x|\apprge A\}} f}{\Leb{\infty}_x} +A^{s-4}\log (A)\norm{f}{\Leb{\infty}_x}\apprle \norm{\action{x}^{2s} f}{\Leb{\infty}_x}.
 \end{align*}

 Therefore, we get the desired claim \eqref{eq:As-estimate} by combining \eqref{eq:low-frequency-1} and \eqref{eq:high-frequency-1}.
 \end{proof}

\bibliographystyle{amsplain}
\bibliography{Refs}

\end{document}